\documentclass{amsart}
\usepackage{amsmath}
\usepackage{amssymb}
\usepackage{amsthm}
\usepackage{graphicx}
\newtheorem{thm}{Theorem}[section]

\newtheorem{cor}[thm]{Corollary}

\theoremstyle{definition}

\newtheorem{note}[thm]{Note}

\usepackage{kotex} 
\begin{document}
	
	\title[The Principle of Indivisible Integrity]{On the Incompatibility of Rearrangement with Convergence: An Axiomatic Approach to Holomorphic Recurrence Relations}
	
	\author{Yoon-Seok Choun}
	\address{Department of Physics, POSTECH, Pohang, Gyeongbuk 37673, South Korea}
	\email{ychoun@gmail.com}
	
	\begin{abstract}  
		In classical analysis, the convergence behavior of power series solutions to differential or recurrence equations is assumed to be invariant under internal rearrangement. This paper challenges that foundational belief by proving that, for holomorphic solutions to higher-order recurrence relations (order $\geq$ 3), rearrangement of internal terms systematically reduces the radius of convergence. This contradicts key assumptions underlying Fuchs’ theorem and the Poincaré--Perron theorem.
		
		To resolve this discrepancy, I propose the \textit{Principle of Indivisible Integrity}, an axiom that restricts arbitrary reordering within analytic computations. The paper provides both analytic proofs and numerical counterexamples (Theorem 3.3, Table 3) showing that the violation of this principle leads to structural divergence, even under conditions classically deemed convergent.
		
		The implications are wide-ranging: this framework calls for a reevaluation of how analytic structures are treated in recurrence-based methods across quantum mechanics, general relativity, and spectral theory. It also raises essential questions about the philosophical foundations of computation and mathematical rigor in the age of machine-optimized symbolic processing.
		
		At its core, this work is not merely a correction to existing theorems, but a philosophical appeal: to recognize that in an era where speed dominates meaning, the integrity of mathematical order must be preserved not by efficiency, but by principle.
	\end{abstract}

	\maketitle
	
	\section{Introduction}  
	
Power series solutions to differential equations have long served as essential tools in mathematical physics, especially in the analysis of singularities and recurrence structures. Classical results such as Fuchs’ theorem and the Poincaré–Perron theorem offer established criteria for determining convergence radii. These foundational theorems, however, carry an implicit assumption: that the internal rearrangement of terms does not affect convergence. While this assumption holds for second-order recurrence relations, it may break down in higher-order cases where the structural dependency among terms becomes nontrivial.

This paper introduces the \textit{Principle of Indivisible Integrity}, which asserts that, in higher-order recurrence relations, terms must be treated not as computational elements but as indivisible structural units—analogous to points in the complex plane. Just as a complex number $z = x + iy$ retains its meaning only when treated as an inseparable entity, the convergence behavior of recurrence-based solutions depends on maintaining the structural integrity of their components.

This principle marks a fundamental departure from classical assumptions:

\begin{enumerate}
	\item While Fuchs’ theorem relates singularities to convergence, the proposed principle highlights how internal structure within recurrence terms can itself govern convergence behavior.
	\item The Poincaré–Perron theorem assumes that rearrangement does not affect convergence. We demonstrate that this assumption fails in higher-order settings, where multiple structural components interact.
	\item Classical theory views recurrence terms as algebraic or numerical inputs. The new principle reframes them as structural primitives, whose integrity is essential for preserving analytic validity.
\end{enumerate}

This reinterpretation not only resolves latent inconsistencies in classical frameworks but also establishes a unified perspective on how structural constraints govern analytic behavior. By recognizing recurrence terms as indivisible, we derive a more coherent and rigorous foundation for analyzing convergence in power series.

We show that while second-order recurrence relations maintain absolute convergence under arbitrary term rearrangement, higher-order relations embed subcomponents whose dislocation alters the radius of convergence. This reveals a previously overlooked limitation in classical convergence theory.

To address this, we introduce structural constraints that define the boundaries of stable convergence. These constraints clarify when classical results remain valid, and when they require refinement. The framework is applied to Heun-type differential equations, revealing how structural integrity affects both analytic form and numerical convergence.

This work addresses two guiding questions:
\begin{enumerate}
	\item How does the Principle of Indivisible Integrity reshape the structural criteria for convergence in higher-order differential equations?
	\item What mathematical conditions are required to preserve convergence behavior under complex recurrence constructions?
\end{enumerate}

By extending the classical framework, this paper offers a structural rethinking of recurrence analysis, revealing latent dependencies and proposing an axiomatic refinement that may reshape future approaches to analytic solutions.   
	 \section{Motivation}
	 
	 The classical framework of power series convergence rests on the assumption that internal rearrangement of terms does not alter the radius of convergence. While this principle holds for second-order recurrence relations—where terms behave as algebraically independent—it fails to generalize to higher-order relations, where the coefficients are constructed from multiple interdependent terms.
	 
	 This paper challenges the generality of this assumption by proposing the \textit{Principle of Indivisible Integrity}, which asserts that in higher-order recurrence relations, the radius of convergence can be distorted if the structural integrity of fundamental terms is violated through decomposition or rearrangement. In such settings, term independence no longer holds; instead, intrinsic structural dependencies emerge between components of the coefficients.
	 
	 By formulating structural constraints on absolute convergence, we provide a framework to analyze how implicit term rearrangement affects the analytic stability of power series. This principle extends classical convergence analysis by identifying previously unrecognized structural limitations that govern the behavior of solutions to higher-order differential equations.
	 
	 The key insight is this: to preserve the analytical validity of power series expansions, recurrence terms must be treated as indivisible structural units. The following sections provide mathematical proofs demonstrating how term rearrangement can reduce the radius of convergence and justify the necessity of the no-resummation constraint.

	\section{Proof of the Breakdown of Absolute Convergence in Higher-Order Recurrence Relations}  
	In second-order recurrence relations, the absolute convergence of power series guarantees that the radius of convergence remains invariant under term rearrangement. However, this invariance does not hold for higher-order recurrence relations due to the internal structure of the recurrence coefficients.
	
	To formalize this distinction, we introduce the concept of an \textit{Indivisible Term}, denoted by $d_i$, which represents a composite expression formed by multiple recurrence coefficients $A_n, B_m, \dots$. Specifically,
	\begin{equation}
		d_i = \sum\left( \prod_{i=0}^{j}\prod_{l=0}^{k}\cdots A_{i}B_{j}\cdots\right).
	\end{equation}
	While absolute convergence typically implies stability under term rearrangement, we demonstrate that rearranging the internal components of $d_i$---a process referred to as \textit{Internal Component Rearrangement}---can cause a contraction in the radius of convergence.

	\textbf{Definition (Breakdown of Absolute Convergence in Higher-Order Recurrence Relations):}
	A power series solution is said to exhibit a breakdown of absolute convergence if internal component rearrangement leads to a reduction in the radius of convergence, despite the recurrence relation being satisfied. This necessitates treating each $d_i$ as an indivisible term to preserve structural integrity.
	
	By analyzing the minimal domain under internal component rearrangement, we establish the necessity of structural constraints to uphold the convergence behavior predicted by classical results such as Fuchs' theorem. These findings highlight the essential distinction between second-order and higher-order recurrence relations in terms of convergence stability. 
	\begin{note}
		For a fixed $k\in \mathbb{N}$, there is the $(k+1)$-term recurrence relation with constant coefficients such as
		\begin{equation}
			c_{k+1,n+1} = \alpha_1 c_{k+1,n} + \alpha_2 c_{k+1,n-1} + \alpha_3 c_{k+1,n-2} + \cdots + \alpha_k c_{k+1,n-k+1}
			\label{eq:19}
		\end{equation}
		with seed values $c_{k+1,j} = \sum_{i=1}^{j}\alpha_i c_{k+1,j-i}$ where $j=1,2,3,\cdots,k-1$ \& $k\in \mathbb{N}-\{1\}$. And $c_{k+1,0}=1$ is chosen for simplicity from now on.
		The generating function of the sequence of (\ref{eq:19}) is given by
		\begin{equation}
			y(x) = \sum_{n=0}^{\infty} c_{k+1,n} x^n = \sum_{r=0}^{\infty} \left(\sum_{j=1}^{k}\alpha_j x^j\right)^r
			\label{eq:21}
		\end{equation}
		The domain of absolute convergence of (\ref{eq:21}) is written by
		\begin{equation}
			\mathcal{D} := \left\{ x \in \mathbb{C} \,\bigg|\, \sum_{j=1}^{k}|\alpha_j x^j| < 1 \right\}
			\label{eq:22}
		\end{equation}
	\end{note}
	\begin{proof}
		See Section 3.3 in Flajolet and Sedgewick~\cite{Sedg1996}. A summation series expansion of (\ref{eq:21}) is
		\begin{equation}
			\sum_{n=0}^{\infty } c_{k+1, n} x^n = \sum_{i_1 =0}^{\infty } \sum_{i_2 =0}^{\infty }\cdots \sum_{i_k =0}^{\infty} \frac{(i_1 +i_2 +\cdots +i_k)!}{i_1 !\;i_2!\; \cdots i_k!}\left( \alpha _{1} x\right)^{i_1 } \left( \alpha _{2} x^2\right)^{i_2 } \cdots \left( \alpha _{k} x^k\right)^{i_k }  \label{cc:1}
		\end{equation}
		A real (or complex) series $\sum_{n=0}^{\infty } u_n$ is referred to converge absolutely if the series of moduli $\sum_{n=0}^{\infty } |u_n|$  converge. And the series of absolute values (\ref{cc:1}) is
		\begin{equation}
			\sum_{i_1 =0}^{\infty } \sum_{i_2 =0}^{\infty }\cdots \sum_{i_k =0}^{\infty} \frac{(i_1 +i_2 +\cdots +i_k)!}{i_1 !\;i_2!\; \cdots i_k!}\left| \alpha _{1} x\right|^{i_1 } \left| \alpha _{2} x^2\right|^{i_2 } \cdots \left| \alpha _{k} x^k\right|^{i_k }  = \sum_{r=0}^{\infty } \left(\sum_{j=1}^{k}|\alpha_j x^j|\right)^r  \nonumber
		\end{equation}
		This $k$-tuple series is absolutely convergent for $ \sum_{j=1}^{k}|\alpha_j x^j|  <1$.
		\qed  
	\end{proof} 
	A Fuchsian differential equation of order $j$ with variable coefficients is of the form
	\begin{equation}
		a_j(x)y^{(j)}+ a_{j-1}(x)y^{(j-1)}+ \cdots + a_1(x)y^{'} + a_0(x)y =0
		\label{jj:1}
	\end{equation}
	Assuming its solution as a Frobenius series in the form
	\begin{equation}
		y(x) = x^{\lambda}\sum_{n=0}^{\infty } d_n x^{n}
		\label{jj:2}
	\end{equation}
	where $d_0 =1$ chosen for simplicity from now on. $\lambda $ is an indicial root. For a fixed $k\in \mathbb{N}$,  we suggest that we obtain the $(k+1)$-term recurrence relation putting (\ref{jj:2}) in (\ref{jj:1}).
	\begin{equation}
		d_{n+1} = \alpha _{1,n} d_{n} + \alpha _{2,n} d_{n-1} + \alpha _{3,n} d_{n-2} + \cdots + \alpha _{k,n} d_{n-k+1}
		\label{jj:3}
	\end{equation}
	with seed values $d_{j} = \sum_{i=1}^{j}\alpha _{i,j-1} d_{j-i}  $ where $j=1,2,3,\cdots, k-1$ \& $k\in \mathbb{N}-\{1\}$.
	\begin{thm}
		To ensure the coefficients stabilize asymptotically, we assume $\lim_{n\rightarrow \infty}\alpha_{l,n}= \alpha_l<\infty$ in (\ref{jj:3}), the minimal domain under internal component rearrangement for (\ref{jj:2}) is given by:
		\begin{equation}
			\mathcal{D} := \left\{ x \in \mathbb{C} \,\bigg|\, \sum_{m=1}^{k}|\alpha_m x^m| < 1 \right\}
			\label{eq:min_domain}
		\end{equation}
		where $\mathcal{D}$ represents the minimal domain under internal component rearrangement.
		\label{main}
	\end{thm}  
	\begin{thm}
		To ensure the coefficients stabilize asymptotically, we assume
		
		${\displaystyle \lim_{n\rightarrow \infty}\alpha _{l,n}= \alpha _{l}<\infty }$ in (\ref{jj:3}),  the minimal domain under internal component rearrangement for (\ref{jj:2}) is given by:
		\begin{equation}
			\mathcal{D} :=  \left\{  x \in \mathbb{C} \Bigg| \sum_{m=1}^{k}\left| \alpha _{m} x^m \right| <1 \right\}
			\label{jj:4}
		\end{equation}\label{main}
		where $\mathcal{D}$ represents the minimal domain under internal component rearrangement.
	\end{thm}
	\begin{proof}
		As shown in the introduction, the concept of absolute Convergence plays a crucial role in determining the minimal domain under internal component rearrangement. This theorem formalizes the relationship between term rearrangement and the radius of convergence.
		Here, absolute convergence refers to the altered convergence domain after term rearrangement. To prove this for all $k \in \mathbb{N}$, consider the $(k+1)$-term recurrence relation where $k = 1, 2, 3, \dots$. Assuming ${\displaystyle \lim_{n\rightarrow \infty}\alpha _{l,n}= \alpha _{l}<\infty }$, let $\alpha _{l,n} = \alpha _{l} \overline{\alpha _{l,n}}$.
		
		Here, $\overline{\alpha _{l,n}}$ represents a normalized coefficient satisfying $|\overline{\alpha _{l,n}}| < 1 + \epsilon$ for $n \geq N$ and any positive error bound $\epsilon$. Let $\tilde{\alpha _{l}} = (1 + \epsilon)\alpha _{l}$ for simplicity. Following the cases for 2-term, 3-term, and 4-term recurrence relations, we generalize the minimal domain under internal component rearrangement condition for a $(k+1)$-term recurrence relation as follows. 
		
		\subsubsection*{Case A: The 2-term recurrence relation with non-constant coefficients of a linear ODE} 
		
		In general, 2-term recurrence relation putting (\ref{jj:2}) in (\ref{jj:1}) is given by
		\begin{equation}
			d_{n+1} = \alpha _{1,n} d_{n}
			\label{jj:5}
		\end{equation}
		If $\left|\overline{\alpha _{1,n}}\right| < 1+\epsilon $ when $n\geq N$, then the series of absolute values, $1+|d_1||x|+|d_2||x|^2+ |d_3||x|^3 +\cdots$, is dominated by the convergent series
		\begin{multline}
			1+|d_1||x|+|d_2||x|^2+ |d_3||x|^3 +\cdots +|d_N||x|^N \left\{ 1+ |\tilde{\alpha _{1}}| |x|+ |\tilde{\alpha _{1}}|^2 |x|^2+ |\tilde{\alpha _{1}}|^3 |x|^3+ \cdots \right\}\\= \sum_{n=0}^{N-1}|d_n||x|^n +\frac{|d_N||x|^N}{1-|\tilde{\alpha _{1}} x|}
			\label{jj:6}
		\end{multline}
		(\ref{jj:6}) is absolutely convergent for $|\tilde{\alpha _{1}} x|<1$. We know $ |\alpha _{1} x|<|\tilde{\alpha _{1}} x|$. Therefore, a series for the 2-term recurrence relation is convergent for  $|\alpha _{1} x|<1$.
		
		\subsubsection*{Case B: The 3-term recurrence relation with non-constant coefficients of a linear ODE} 
		
		The 3-term recurrence relation putting (\ref{jj:2}) in (\ref{jj:1}) is given by
		\begin{equation}
			d_{n+1}= \alpha _{1,n} \;d_n + \alpha _{2,n} \;d_{n-1} \hspace{1cm};n\geq 1
			\label{jj:7}
		\end{equation}
		with seed values $d_1= \alpha _{1,0} $ and $d_0=1$.
		When $n\geq N$, we have $\left|\overline{\alpha _{1,n}}\right|, \left|\overline{\alpha _{2,n}}\right| < 1+\epsilon $.
		
		For $n=N, N+1, N+2, \cdots$ in succession, take the modulus of the general term of $d_{n+1}$ in  (\ref{jj:7}) 
		\begin{multline}
			|d_{N+1+j}| = |\alpha _{1,N+j}| |d_{N+j}| + |\alpha _{2,N+j}|  |d_{N-1+j}|  \leq  |\tilde{\alpha _{1}}| |d_{N+j}| + |\tilde{\alpha _{2}}|  |d_{N-1+j}| \\
			\leq c_{3,j+1}  |d_{N}| + c_{3,j} |\tilde{\alpha _{2}}| |d_{N-1}|
			\label{jj:8}
		\end{multline}
		Here, $c_{3,j}$ where $j \in  \mathbb{N}_0$ is a generated coefficient s of the 3-term recurrence relation  in  (\ref{eq:19}) with taking $\alpha _{1} \rightarrow |\tilde{\alpha _{1}}| $ and $\alpha _{2} \rightarrow |\tilde{\alpha _{2}}| $.
		
		According to (\ref{jj:8}), then the series of absolute values, $1+|d_1||x|+|d_2||x|^2+ |d_3||x|^3 +\cdots$, is dominated by the convergent series
		\begin{multline}
			\sum_{n=0}^{N-1}|d_n||x|^n + |d_N||x|^N + \sum_{k=0}^{\infty }\left( c_{3,k+1} |d_N| + c_{3,k}|\tilde{\alpha _{2}}| |d_{N-1}|  \right) |x|^{N+1+k}  \\
			= \sum_{n=0}^{N-1}|d_n||x|^n + \left( |d_N| + |\tilde{\alpha _{2}}| |d_{N-1}| |x|\right) |x|^N \sum_{i=0}^{\infty }c_{3,i} |x|^{i}
			\label{jj:9}
		\end{multline}
		From the generating function for the 3-term recurrence relation in (\ref{eq:21}),  we conclude
		\begin{equation}
			\sum_{i=0}^{\infty }c_{3,i} |x|^{i} =  
			\sum_{r=0}^{\infty } \left(|\tilde{\alpha _{1}}x| + |\tilde{\alpha _{2}}x^2| \right)^r 
			\label{jj:10}
		\end{equation}
		(\ref{jj:10}) is absolute convergent for $|\tilde{\alpha _{1}} x|+ |\tilde{\alpha _{2}}x^2|<1$. We know $ |\alpha _{l} x^l|<|\tilde{\alpha _{l}} x^l|$ where $l=1,2$ . Therefore, a series for the 3-term recurrence relation is convergent for  $|\alpha _{1} x|+|\alpha _{2} x^2| <1$.
		
		\subsubsection*{Case C: The 4-term recurrence relation with non-constant coefficients of a linear ODE} 
		
		The 4-term recurrence relation putting (\ref{jj:2}) in (\ref{jj:1}) is given by
		\begin{equation}
			d_{n+1}= \alpha _{1,n} \;d_n + \alpha _{2,n} \;d_{n-1} + \alpha _{3,n} \;d_{n-2} \hspace{1cm};n\geq 2
			\label{jj:11}
		\end{equation}
		with seed values $d_1= \alpha _{1,0} $, $d_2= \alpha _{1,0}\alpha _{1,1}+ \alpha _{2,1} $ and $d_0=1$.
		When $n\geq N$, we have $\left|\overline{\alpha _{1,n}}\right|, \left|\overline{\alpha _{2,n}}\right|, \left|\overline{\alpha _{3,n}}\right| < 1+\epsilon $.
		
		For $n=N, N+1, N+2, \cdots$ in succession, take the modulus of the general term of $d_{n+1}$ in  (\ref{jj:11}) 
		\begin{multline}
			|d_{N+1+j}| = |\alpha _{1,N+j}| |d_{N+j}| + |\alpha _{2,N+j}|  |d_{N-1+j}| + |\alpha _{3,N+j}|  |d_{N-2+j}|  \\ 
			\leq  |\tilde{\alpha _{1}}| |d_{N+j}| + |\tilde{\alpha _{2}}|  |d_{N-1+j}|  + |\tilde{\alpha _{3}}|  |d_{N-2+j}|    \\
			\leq c_{4,j+1}  |d_{N}| +\left( c_{4,j+2}- c_{4,j+1} |\tilde{\alpha _{1}}| \right)   |d_{N-1}| + c_{4,j} |\tilde{\alpha _{3}}| |d_{N-2}|
			\label{jj:12}
		\end{multline}
		Here, $c_{4,j}$ where $j \in  \mathbb{N}_0$ is a  generated coefficient  of the 4-term recurrence relation in  (\ref{eq:19}) with taking $\alpha _{l} \rightarrow |\tilde{\alpha _{l}}| $ where $l=1,2,3$.
		
		According to (\ref{jj:12}), then the series of absolute values, $1+|d_1||x|+|d_2||x|^2+ |d_3||x|^3 +\cdots$, is dominated by the convergent series
		\begin{multline}
			\sum_{n=0}^{N-1}|d_n||x|^n + |d_N||x|^N + \sum_{k=0}^{\infty }\bigg( c_{4,k+1} |d_N| + \left( c_{4,k+2}-  c_{4,k+1}|\tilde{\alpha _{1}}|  \right) |d_{N-1}| + c_{4,k} |\tilde{\alpha _{3}}||d_{N-2}| \bigg) |x|^{N+1+k}\\
			\leq \sum_{n=0}^{N-1}|d_n||x|^n + |d_N||x|^N + \sum_{k=0}^{\infty }\left( c_{4,k+1} |d_N| +  c_{4,k+2}  |d_{N-1}| + c_{4,k} |\tilde{\alpha _{3}}||d_{N-2}| \right) |x|^{N+1+k} \\
			\leq  \sum_{n=0}^{N-1}|d_n||x|^n + \left( |d_N| + |\tilde{\alpha _{3}}| |d_{N-2}| |x|\right) |x|^N \sum_{i=0}^{\infty }c_{4,i} |x|^{i} + |d_{N-1}||x|^{N-1}\sum_{i=2}^{\infty }c_{4,i} |x|^{i}
			\label{jj:13}
		\end{multline}
		From the generating function for the 4-term recurrence relation in (\ref{eq:21}),  we conclude
		\begin{equation}
			\sum_{i=0}^{\infty }c_{4,i} |x|^{i} =  
			\sum_{r=0}^{\infty } \left(|\tilde{\alpha _{1}}x| + |\tilde{\alpha _{2}}x^2|+ |\tilde{\alpha _{3}}x^3|  \right)^r   \label{jj:14}
		\end{equation}
		(\ref{jj:14}) is absolute convergent for $|\tilde{\alpha _{1}} x|+ |\tilde{\alpha _{2}}x^2| + |\tilde{\alpha _{3}}x^3|<1$. We know $ |\alpha _{l} x^l|<|\tilde{\alpha _{l}} x^l|$ where $l=1,2,3$ . Therefore, a series for the 4-term recurrence relation is convergent for
		$|\alpha _{1} x|+|\alpha _{2} x^2|+|\alpha _{3} x^3|  <1$.
		
		\subsubsection*{Case D: The 5-term recurrence relation with non-constant coefficients of a linear ODE} 
		
		The 5-term recurrence relation putting (\ref{jj:2}) in (\ref{jj:1}) is given by
		\begin{equation}
			d_{n+1}= \alpha _{1,n} \;d_n + \alpha _{2,n} \;d_{n-1} + \alpha _{3,n} \;d_{n-2} + \alpha _{4,n} \;d_{n-3}\hspace{1cm};n\geq 3
			\label{jj:15}
		\end{equation}
		with seed values $d_1= \alpha _{1,0} $, $d_2= \alpha _{1,0}\alpha _{1,1}+ \alpha _{2,1} $, $d_3= \alpha _{1,0}\alpha _{1,1}\alpha _{1,2}+ \alpha _{1,0}\alpha _{2,2}+ \alpha _{1,2} \alpha _{2,1} +\alpha _{3,2} $ and $d_0=1$.
		When $n\geq N$, we have $\left|\overline{\alpha _{1,n}}\right|, \left|\overline{\alpha _{2,n}}\right|, \left|\overline{\alpha _{3,n}}\right|, \left|\overline{\alpha _{4,n}}\right|< 1+\epsilon $.
		
		For $n=N, N+1, N+2, \cdots$ in succession, take the modulus of the general term of $d_{n+1}$ in  (\ref{jj:15}) 
		\begin{multline}
			|d_{N+1+j}| = |\alpha _{1,N+j}| |d_{N+j}| + |\alpha _{2,N+j}|  |d_{N-1+j}| + |\alpha _{3,N+j}|  |d_{N-2+j}| + |\alpha _{4,N+j}|  |d_{N-3+j}| \\
			\leq  |\tilde{\alpha _{1}}| |d_{N+j}| + |\tilde{\alpha _{2}}|  |d_{N-1+j}|  + |\tilde{\alpha _{3}}|  |d_{N-2+j}|+ |\tilde{\alpha _{4}}|  |d_{N-3+j}| \\
			\leq  c_{5,j+1}  |d_{N}| +\left( c_{5,j+2}- c_{5,j+1} |\tilde{\alpha _{1}}| \right)   |d_{N-1}| + \left( c_{5,j+3}-
			\left( c_{5,j+2} |\tilde{\alpha _{1}}|+ c_{5,j+1} |\tilde{\alpha _{2}}|\right) \right)  |d_{N-2}| + c_{5,j} |\tilde{\alpha _{4}}| |d_{N-3}|
			\label{jj:16}
		\end{multline}
		Here, $c_{5,j}$ where $j \in  \mathbb{N}_0$ is a  generated coefficient of the 5-term recurrence relation in  (\ref{eq:19}) with taking $\alpha _{l} \rightarrow |\tilde{\alpha _{l}}| $ where $l=1,2,3,4$.
		
		According to (\ref{jj:16}), then the series of absolute values, $1+|d_1||x|+|d_2||x|^2+ |d_3||x|^3 +\cdots$, is dominated by the convergent series 
		\begin{eqnarray}
			&&\sum_{n=0}^{N-1}|d_n||x|^n + |d_N||x|^N 	   \nonumber\\
			&&+ \sum_{k=0}^{\infty }\bigg( c_{5,k+1} |d_N| + \left( c_{5,k+2}- c_{5,k+1}|\tilde{\alpha _{1}}|   \right) |d_{N-1}| +  \left( c_{5,k+3}- \left(c_{5,k+2} |\tilde{\alpha _{1}}|   +  c_{5,k+1}|\tilde{\alpha _{2}}| \right) \right)|d_{N-2}|\nonumber\\
			&& + c_{5,k}|\tilde{\alpha _{4}}||d_{N-3}|\bigg) |x|^{N+1+k}\nonumber\\
			&&\leq \sum_{n=0}^{N-1}|d_n||x|^n + |d_N||x|^N \nonumber\\
			&&+ \sum_{k=0}^{\infty }\bigg( c_{5,k+1} |d_N| +   c_{5,k+2}  |d_{N-1}| +    c_{5,k+3} |d_{N-2}| + c_{5,k}|\tilde{\alpha _{4}}||d_{N-3}|\bigg) |x|^{N+1+k}  \nonumber\\
			&&\leq  \sum_{n=0}^{N-1}|d_n||x|^n + \left( |d_N| + |\tilde{\alpha _{4}}| |d_{N-3}| |x|\right) |x|^N \sum_{i=0}^{\infty }c_{5,i} |x|^{i} + |d_{N-1}||x|^{N-1} \sum_{i=2}^{\infty }c_{5,i} |x|^{i}  \nonumber\\
			&&+ |d_{N-2}||x|^{N-2}\sum_{i=3}^{\infty } c_{5,i} |x|^{i}  \nonumber
		\end{eqnarray}
		From the generating function for the 5-term recurrence relation in (\ref{eq:21}),  we conclude
		\begin{equation}
			\sum_{i=0}^{\infty }c_{5,i} |x|^{i}  =  
			\sum_{r=0}^{\infty } \left(|\tilde{\alpha _{1}}x| + |\tilde{\alpha _{2}}x^2|+ |\tilde{\alpha _{3}}x^3| + |\tilde{\alpha _{4}}x^4|  \right)^r  
			\label{jj:18}
		\end{equation}
		(\ref{jj:18}) is absolute convergent for $|\tilde{\alpha _{1}} x|+ |\tilde{\alpha _{2}}x^2| + |\tilde{\alpha _{3}}x^3| + |\tilde{\alpha _{4}}x^4|<1$. We know $ |\alpha _{l} x^l|<|\tilde{\alpha _{l}} x^l|$ where $l=1,2,3,4$ . Therefore, a series for the 5-term recurrence relation is convergent for
		$|\alpha _{1} x|+|\alpha _{2} x^2|+|\alpha _{3} x^3| +|\alpha _{4} x^4| <1$.
		
		By the principle of mathematical induction (by repeating similar process for the previous cases such as two, three, four and five term recurrence relations), the series of absolute values for the $(k+1)$-term recursive relation of a linear ODE where $k\in \{ 3,4,5,\cdots\}$, $1+|d_1||x|+|d_2||x|^2+ |d_3||x|^3 +\cdots$, is dominated by the convergent series
		\begin{equation}
			\sum_{n=0}^{N-1}|d_n||x|^n + \bigg( |d_N| + |\tilde{\alpha _{k}}| |d_{N-k+1}| |x|\bigg) |x|^N \sum_{j=0}^{\infty }c_{k+1,j} |x|^{j}
			+ \sum_{i=0}^{k-3}|d_{N-1-i}||x|^{N-1-i}\sum_{j=i+2}^{\infty }c_{k+1,j} |x|^{j}
			\label{jj:19}
		\end{equation}
		Here, $c_{k+1,j}$ where $j \in  \mathbb{N}_0$ is a  generated coefficient of the $(k+1)$-term recurrence relation in  (\ref{eq:19}) with taking $\alpha _{l} \rightarrow |\tilde{\alpha _{l}}| $ where $l=1,2,3,\cdots,k$.
		From the generating function for the $(k+1)$-term recurrence relation in (\ref{eq:21}),  we conclude
		\begin{equation}
			\sum_{j=0}^{\infty }c_{k+1,j} |x|^{j} = \sum_{r=0}^{\infty } \left(\sum_{m=1}^{k}|\tilde{\alpha _{m}} x^m|\right)^r  
			\label{jj:20}
		\end{equation}
		(\ref{jj:20}) is minimal absolute convergence domain under internal component rearrangement for $\left(\sum_{m=1}^{k} |\tilde{\alpha _{m}}x^m|\right) < 1$. We know $ |\alpha _{l} x^l|<|\tilde{\alpha _{l}} x^l|$ where $l=1,2,3, \cdots, k$. Therefore, a series for the $(k+1)$-term recurrence relation of a Fuchsian differential equation is convergent for $\sum_{m=1}^{k} |\alpha _{m}x^m| <1$.
		\qed
	\end{proof} 
	This result offers a rigorous foundation for understanding the intrinsic structural limitations of Fuchs' theorem in the context of power series analysis. It underscores the importance of the \textit{Indivisible Integrity Principle}, which demands preservation of the internal structure of recurrence terms to maintain a stable radius of convergence. By forbidding internal component rearrangement, this principle aligns the behavior of the series with classical theoretical expectations, even in the presence of complex recurrence relations or perturbative constructs.
	
	The \textit{Indivisible Integrity Principle} bridges a crucial gap between classical analysis and modern computational practices, where structural changes can inadvertently impact convergence. Through this formalism, we broaden the theoretical scope of Fuchs' theorem, enabling consistent treatment of power series expansions across analytical, asymptotic, and numerical domains. Beyond addressing a key limitation in traditional frameworks, this work provides a mathematically robust methodology for preserving structural coherence in power series solutions within diverse mathematical fields, including differential equations, function theory, and applied analysis.

	\section{Application to the Heun's Differential Equation}\label{sec.1}
	
	The Heun equation, known as the most general Fuchsian equation with four regular singularities, arises in various domains of modern physics, including quantum gravity, general relativity, and molecular spectroscopy. Unlike hypergeometric functions which involve two-term recurrence relations, the solutions to the Heun equation often involve more complex three-term recurrence relations. This structural complexity necessitates a more refined approach to analyzing convergence behavior.
	
	Traditionally, the Poincar\'e--Perron (P--P) theorem has been used to determine the radius of convergence for power series solutions to Fuchsian equations, including the Heun equation. However, our study demonstrates that the radius derived from the P--P theorem only accounts for conditional convergence and fails to address absolute convergence. In particular, numerical results consistently show that the radius obtained via Theorem~\ref{main} is smaller than that obtained using the P--P theorem, due to the latter's implicit assumption of non-rearranged coefficients.
	
	By leveraging Theorem~\ref{main}, we establish the radius of absolute convergence for the power series solution to the Heun equation and highlight the significant distinctions in the convergence domains. Furthermore, this framework reveals why the P--P theorem overlooks absolute convergence by examining its foundational assumptions. Consequently, this approach provides a more precise framework for interpreting the behavior of Fuchsian equations with multi-term recurrence relations, especially in modern physical contexts that involve higher-order differential equations.
	
	The Heun equation is a second-order linear ODE expressed as \cite{Heun1889}:
	\begin{equation}
		\frac{d^2y}{dx^2} + \left(\frac{\gamma}{x} + \frac{\delta}{x-1} + \frac{\varepsilon}{x-a}\right)\frac{dy}{dx} + \frac{\alpha\beta x-q}{x(x-1)(x-a)}y = 0
		\label{eq:1}
	\end{equation}
	with the constraint \(\varepsilon = \alpha + \beta - \gamma - \delta + 1\). Here, the parameter \(a \neq 0\) denotes the singularity, \(\alpha, \beta, \gamma, \delta, \varepsilon\) are exponent parameters, and \(q\) is the accessory parameter. \(\alpha\) and \(\beta\) are symmetric, giving a total of six free parameters. The equation features four regular singular points located at \(0\), \(1\), \(a\), and \(\infty\), with corresponding exponents \(\{0, 1-\gamma\}\), \(\{0, 1-\delta\}\), \(\{0, 1-\varepsilon\}\), and \(\{\alpha, \beta\}\).
	
	To study convergence, we adopt the restriction \(|a| > 1\) for the local Heun function near \(x = 0\). Nevertheless, we allow \(a \in \mathbb{R}\) or \(a \in \mathbb{C}\) to examine how the convergence radius behaves across different parameter settings. Assuming a power series solution:
	\begin{equation}
		y(x) = x^{\lambda}\sum_{n=0}^{\infty}d_n x^n
		\label{eq:2}
	\end{equation}
	and substituting it into \eqref{eq:1} yields the recurrence:
	\begin{equation}
		d_{n+1} = \alpha_{1,n}d_n + \alpha_{2,n}d_{n-1} \quad;\,n\geq 1
		\label{eq:3}
	\end{equation}
	where \(\alpha_{1,n} = \alpha_1\overline{\alpha_{1,n}}\), \(\alpha_{2,n} = \alpha_2\overline{\alpha_{2,n}}\), and the specific parameter forms are given in \eqref{eq:4}. We also identify two indicial roots: \(\lambda = 0\) and \(\lambda = 1-\gamma\).

		\begin{align}
		\alpha_1 &= \frac{1+a}{a} \nonumber \\
		\alpha_2 &= -\frac{1}{a} \nonumber \\
		\overline{\alpha_{1,n}} &= \frac{n^2 + \frac{\gamma + \varepsilon - 1 + 2\lambda + a(\gamma + \delta - 1 + 2\lambda)}{1+a}n + \frac{\lambda(\gamma + \varepsilon - 1 + \lambda + a(\gamma + \delta - 1 + \lambda)) + q}{1+a}}{n^2 + (1 + \gamma + 2\lambda)n + (1 + \lambda)(\gamma + \lambda)} \nonumber \\
		&= \frac{n^2 + \frac{\alpha + \beta - \delta + 2\lambda + a(\gamma + \delta - 1 + 2\lambda)}{1+a}n + \frac{\lambda(\alpha + \beta - \delta + \lambda + a(\gamma + \delta - 1 + \lambda)) + q}{1+a}}{n^2 + (1 + \gamma + 2\lambda)n + (1 + \lambda)(\gamma + \lambda)} \nonumber \\
		\overline{\alpha_{2,n}} &= \frac{n^2 + (\alpha + \beta - 2 + 2\lambda)n + (\alpha - 1 + \lambda)(\beta - 1 + \lambda)}{n^2 + (\gamma + 1 + 2\lambda)n + (1 + \lambda)(\gamma + \lambda)} \nonumber \\
		d_1 &= \overline{\alpha_{1,0}}d_0 = \alpha_1\overline{\alpha_{1,0}}d_0
		\label{eq:4}
	\end{align}

	The following subsections explore the classical framework provided by the P--P theorem and contrast it with our absolute convergence analysis using internal component rearrangement, culminating in a deeper understanding of structural stability in Heun-type solutions. 
	\subsection{P--P Theorem and Its Applications for Solutions of Power Series}\label{sec.2}
	
	This section reviews the asymptotic behavior of solutions to linear difference equations with constant coefficients. Consider a linear recurrence relation of order $k+1$ with constant coefficients $\alpha_i$, where $i=0,1,2,\dots,k$:
	\begin{equation}
		u(n+1) + \alpha_1 u(n) + \alpha_2 u(n-1) + \alpha_3 u(n-2) + \cdots + \alpha_k u(n-k+1) = 0
		\label{eq:pp1}
	\end{equation}
	where $\alpha_k \neq 0$.
	
	The characteristic polynomial of the recurrence relation \eqref{eq:pp1} is given by:
	\begin{equation}
		\rho^k + \alpha_1\rho^{k-1} + \alpha_2\rho^{k-2} + \cdots + \alpha_k = 0
		\label{eq:pp2}
	\end{equation}
	The roots of the characteristic equation \eqref{eq:pp2} are denoted by $\lambda_1,\dots,\lambda_k$.
	
	In 1885, H. Poincaré stated that
	\begin{equation}
		\lim_{n\rightarrow \infty} \frac{u(n+1)}{u(n)} \nonumber
	\end{equation}
	is equal to one of the roots of the characteristic equation \cite{Poin1885}. This result was extended by O. Perron in 1921 \cite{Perr1921}:
	The Poincaré–Perron theorem assumes conditional convergence and does not guarantee absolute convergence, which is necessary to preserve structural invariance and uniqueness.
	\begin{thm}[P--P Theorem \cite{Miln1933}]\label{thm.2}
		If the coefficient of $u(n)$ in a difference equation of order $k$ is nonzero for $n=0,1,2,\dots$, then the equation has $k$ fundamental solutions $u_1(n),\dots,u_k(n)$ such that:
		\begin{equation}
			\lim_{n\rightarrow \infty} \frac{u_i(n+1)}{u_i(n)} = \lambda_i
		\end{equation}
		where $i=1,2,\dots,k$, $\lambda_i$ are the roots of the characteristic equation, and $n\rightarrow \infty$ in positive integer increments.
	\end{thm}
	
	A recurrence relation for the coefficients emerges when substituting a series $y(x)=\sum_{n=0}^\infty d_n x^n$ into a Fuchsian equation. In general, a three-term recurrence relation is expressed as:
	\begin{equation}
		d_{n+1} = \alpha_{1,n}d_n + \alpha_{2,n}d_{n-1} \quad;\,n \geq 1
		\label{eq:pp3}
	\end{equation}
	with initial values $d_1= \alpha_{1,0}d_0$. For the asymptotic behavior of \eqref{eq:pp3}, $\lim_{n\rightarrow \infty} \alpha_{j,n} = \alpha_j < \infty$ $(j=1,2)$ exists.
	
	The asymptotic recurrence relation is then:
	\begin{equation}
		\overline{d}_{n+1} = \alpha_1\overline{d}_n + \alpha_2\overline{d}_{n-1} \quad;\,n \geq 1
		\label{eq:pp4}
	\end{equation}
	where $\overline{d}_1= \alpha_1\overline{d}_0$ and $\overline{d}_0 = 1$. From the P--P theorem, the characteristic polynomial of the recurrence relation is:
	\begin{equation}
		\rho^2 - \alpha_1\rho - \alpha_2 = 0
		\label{eq:pp5}
	\end{equation}
	
	The roots of the polynomial \eqref{eq:pp5} have two moduli:
	\begin{equation}
		\rho_1 = \frac{\alpha_1 - \sqrt{\alpha_1^2 + 4\alpha_2}}{2}, \quad
		\rho_2 = \frac{\alpha_1 + \sqrt{\alpha_1^2 + 4\alpha_2}}{2}
		\label{eq:root}
	\end{equation}
	Thus, $\lim_{n\rightarrow \infty} |d_{n+1}/d_n| = \lim_{n\rightarrow \infty} |\overline{d}_{n+1}/\overline{d}_n|$.
	
	In general:
	\begin{itemize}
		\item If $|\rho_1| < |\rho_2|$, then $\lim_{n\rightarrow \infty} |\overline{d}_{n+1}/\overline{d}_n| = |\rho_2|$, and the radius of convergence is $|\rho_2|^{-1}$.
		\item If $|\rho_2| < |\rho_1|$, then $\lim_{n\rightarrow \infty} |\overline{d}_{n+1}/\overline{d}_n| = |\rho_1|$, and the radius of convergence is $|\rho_1|^{-1}$.
	\end{itemize}
	
	In special cases:
	\begin{itemize}
		\item If $|\rho_1| = |\rho_2|$ and $\rho_1 \neq \rho_2$, the series diverges.
		\item If $\rho_1 = \rho_2$, the series converges.
	\end{itemize}
	
	More details are provided in Appendix B of Part A \cite{Ronv1995}, Wimp (1984) \cite{Wimp1984}, Kristensson (2010) \cite{Kris2010}, and Erdélyi (1955) \cite{Erde1955}.
	
	Table~\ref{cb.1} lists all possible boundary conditions of a local Heun function at $x=0$, corresponding to the exponent zero, where $a,x \in \mathbb{R}$. These are obtained by applying the P--P theorem with $\alpha_1 = \frac{1+a}{a}$ and $\alpha_2 = -\frac{1}{a}$ in \eqref{eq:root}.
	
	\begin{table}[h]
		\centering
		\begin{tabular}{|l|c|}
			\hline
			Range of the coefficient $a$ & Range of the independent variable $x$ \\
			\hline
			$a = 0$ & No solution \\
			$|a| \geq 1$ & $|x| < 1$ \\
			$|a| < 1$ & $|x| < |a|$ \\
			\hline
		\end{tabular}
		\caption{Boundary conditions of $x$ for a local Heun function at $x=0$ using the P--P theorem.}
		\label{cb.1}
	\end{table}
	\begin{figure}[h]
		\centering
		\includegraphics[scale=.35]{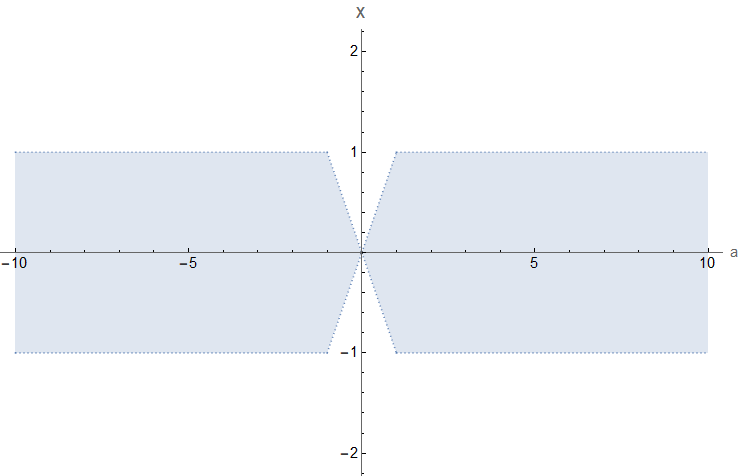}
		\caption{Domain of convergence according to the original P--P theorem.}
		\label{choun2}
	\end{figure}
	Figure~\ref{choun2} represents the convergence domain of the series for the Heun equation at $x=0$. The shaded region indicates convergence, excluding the dotted lines.
	
	Finally, the power series solution derived using the P--P theorem is expressed as:
	\begin{align}
		y^P(x) &= \sum_{n=0}^{\infty}d_n x^{\lambda} \nonumber \\
		&= x^{\lambda}\Big(1 + \alpha_{1,0}x + (\alpha_{2,1} + \alpha_{1,0}\alpha_{1,1})x^2 + \cdots\Big)
		\label{eq:ppsol}
	\end{align}
	\subsection{Structural Modifications and Radius of Convergence in a Power Series with a Three-Term Recurrence Relation}\label{sec.7}
	
	By applying Thm.~\ref{main}, the condition of absolute convergence for a local Heun function at $x = 0$, where $\alpha_1 = (1 + a)/a$ and $\alpha_2 = -1/a$, is given by:
	\begin{equation}
		\left|\frac{1+a}{a}x\right| + \left|-\frac{1}{a}x^2\right| < 1
		\label{eq:17}
	\end{equation}
	
	The coefficient $a$ determines the range of the independent variable $x$ as shown in \eqref{eq:17}. The precise ranges of $a$ and $x$ are given in Table~\ref{cb.2}, where $a,x \in \mathbb{R}$.
	\begin{table}[htbp]
		\centering
		\begin{tabular}{|l|c|}
			\hline
			Range of the coefficient $a$ & Range of the independent variable $x$ \\
			\hline
			$a = 0$ & No solution \\
			$a > 0$ & $|x| < \frac{1}{2}(-1-a+\sqrt{a^2+6a+1})$ \\
			$-1 < a < 0$ & $a < x < -a$ \\
			$a \leq -1$ & $|x| < 1$ \\
			\hline
		\end{tabular}
		\caption{Boundary conditions of $x$ for the infinite series of a local Heun function at $x=0$.}
		\label{cb.2}
	\end{table}
	\begin{figure}[h]
		\centering
		\includegraphics[scale=.35]{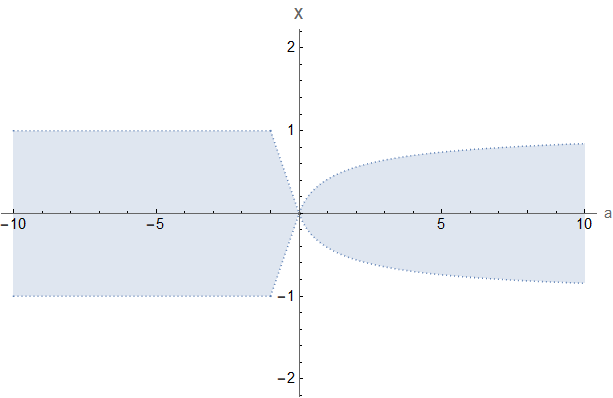}
		\caption{Structural Constraints in the P--P Theorem.}
		\label{choun1}
	\end{figure}
	The radius of convergence from Table~\ref{cb.2} is depicted as the shaded region in Fig.~\ref{choun1}. It does not include the dotted lines, no solution exists at the origin (black point), and the maximum modulus of $x$ is unity.
	
	The dotted boundary lines for the shaded area where $a > 0$ in Fig.~\ref{choun1} are determined by:
	\begin{equation}
		\lim_{a\rightarrow N} \frac{-1-a+\sqrt{a^2+6a+1}}{2} \sim 1,
	\end{equation}
	where $N$ is a sufficiently large positive real number. Hence, it can be argued that $|x| < 1$ as $a \to N$. For example, if $a = 10$, then $|x| < 0.84429$, and if $a = 100$, then $|x| < 0.98058$.
	
	By rearranging the coefficients $\alpha_{1,n}$ and $\alpha_{2,n}$ in each sequence $d_n$ in \eqref{eq:3}, where $d_0 = 1$ for simplicity, a local Heun series solution can be expressed as:
	\begin{equation}
		y^A(x) = x^{\lambda}\left(y_0^A(x) + y_1^A(x)\eta + \sum_{\tau=2}^{\infty}y_{\tau}^A(x)\eta^{\tau}\right),
		\label{eq:heun_series}
	\end{equation}
	where
	\begin{align*}
		y_0^A(x) &= \sum_{i_0=0}^{\infty}\prod_{i_1=0}^{i_0-1}\overline{\alpha_{2,2i_1+1}}z^{i_0}, \\
		y_1^A(x) &= \sum_{i_0=0}^{\infty}\overline{\alpha_{1,2i_0}}\prod_{i_1=0}^{i_0-1}\overline{\alpha_{2,2i_1+1}}
		\sum_{i_2=i_0}^{\infty}\prod_{i_3=i_0}^{i_2-1}\overline{\alpha_{2,2i_3+2}}z^{i_2}, \\
		y_{\tau}^A(x) &= \sum_{\tau=2}^{\infty}\Bigg\{\sum_{i_0=0}^{\infty}\overline{\alpha_{1,2i_0}}
		\prod_{i_1=0}^{i_0-1}\overline{\alpha_{2,2i_1+1}}
		\prod_{k=1}^{\tau-1}\left(\sum_{i_{2k}=i_{2(k-1)}}^{\infty}\overline{\alpha_{1,2i_{2k}+k}}
		\prod_{i_{2k+1}=i_{2(k-1)}}^{i_{2k}-1}\overline{\alpha_{2,2i_{2k+1}+k+1}}\right) \\
		&\quad\times\sum_{i_{2\tau}=i_{2(\tau-1)}}^{\infty}\prod_{i_{2\tau+1}=i_{2(\tau-1)}}^{i_{2\tau}-1}
		\overline{\alpha_{2,2i_{2\tau+1}+\tau+1}}z^{i_{2\tau}}\Bigg\}.
	\end{align*}
	
	Additionally, the variables $\eta$ and $z$ are defined as:
	\begin{equation}
		\eta = \frac{1+a}{a}x, \quad z = -\frac{1}{a}x^2.
	\end{equation}
	The sequence $d_n$ is a combination of $\alpha_{1,n}$ and $\alpha_{2,n}$ terms in \eqref{eq:3}. The series in \eqref{eq:heun_series} is constructed by letting $\alpha_{1,n}$ in $d_n$ be the leading term, observing terms in $d_n$ that include:
	\begin{itemize}
		\item Zero terms of $\alpha_{1,n}$ for the sub-power series $y_0^A(x)$,
		\item One term of $\alpha_{1,n}$ for $y_1^A(x)$,
		\item Two terms of $\alpha_{1,n}$ for $y_2^A(x)$,
		\item Three terms of $\alpha_{1,n}$ for $y_3^A(x)$, and so on.
	\end{itemize}
	
	\subsection{Numerical Comparison of Absolute Convergence and Internal Component Rearrangement}\label{sec.9}
	
	When comparing Table~\ref{cb.2} with Table~\ref{cb.1}, the boundary conditions for the radius of convergence are equivalent for $a < 0$. However, for $a \geq 1$, their ranges of $x$ differ significantly:
	
	\begin{enumerate}
		\item The radius of convergence is unity in Table~\ref{cb.1} at $a=1$, while Table~\ref{cb.2} suggests it is approximately 0.414214.
		
		\item In the region $0 < a < 1$, the maximum absolute value of $x$ differs significantly between Tables~\ref{cb.1} and \ref{cb.2}. For positive $x$, Table~\ref{cb.2} describes a square root function of $a$, with a slope ranging between 0.207107 and 1. In contrast, Table~\ref{cb.1} provides a linear increase with slope 1. For negative $x$, the square root function in Table~\ref{cb.2} has a slope between $-1$ and $-0.207107$, whereas Table~\ref{cb.1} specifies a slope of $-1$.
		
		\item For large $a$, the square root function in Table~\ref{cb.2} approaches $\pm 1$, validating the absolute Convergence theorem. However, for $0 < a < 1$, the absolute Convergence theorem no longer applies to constructing the radius of convergence for local Heun functions.
	\end{enumerate}
	
	The differences between Tables~\ref{cb.1} and \ref{cb.2} can also be analyzed numerically. A sequence $d_n$ is derived by inserting a power series into Heun's equation. The boundary condition of $x$ in Table~\ref{cb.1} is obtained using the ratio $\lim_{n\rightarrow \infty}|d_{n+1}/d_n| = \lim_{n\rightarrow \infty}|\overline{d}_{n+1}/\overline{d}_n|$ in Eq.~\ref{eq:pp4}, while the radius of convergence in Table~\ref{cb.2} is constructed using Theorem~\ref{main}.
	For numerical computations, we set $\overline{\alpha_{1,n}}$ and $\overline{\alpha_{2,n}}$ to unity. The asymptotic recurrence relation for a local Heun function is given by:
	\begin{equation}
		\overline{d}_{n+1} = \alpha_1\overline{d}_n + \alpha_2\overline{d}_{n-1} = \frac{1+a}{a}\overline{d}_n - \frac{1}{a}\overline{d}_{n-1}.
		\label{eq:100}
	\end{equation}
	
	A series solution based on the absolute Convergence theorem is:
	\begin{equation}
		y(x) = \lim_{N\rightarrow\infty}\sum_{n=0}^N \overline{d}_n x^n.
		\label{eq:series1}
	\end{equation}
	
	If Eq.~\eqref{eq:series1} converges absolutely, it can also be expressed as:
	\begin{equation}
		y(x) = \lim_{N\rightarrow\infty}\sum_{n=0}^N \sum_{m=0}^N \frac{(n+m)!}{n!m!}\left(\frac{1+a}{a}x\right)^n\left(-\frac{1}{a}x^2\right)^m.
		\label{eq:series2}
	\end{equation}
	
	For $a=0.8$, the boundary condition in Table~\ref{cb.2} is approximately $-0.368858 < x < 0.368858$, while the radius of convergence in Table~\ref{cb.1} is $-0.8 < x < 0.8$. Substituting $a=0.8$ and $x=0.3$ or $x=0.7$ into Eq.~\eqref{eq:series1} and Eq.~\eqref{eq:series2}, the numerical results are presented in Table~\ref{tab:merged}.
	
	\begin{table}[htbp]
		\centering
		\caption{Numerical Results for $y(x)$ in Different Cases}
		\label{tab:merged}
		\resizebox{\textwidth}{!}{
			\begin{tabular}{c|c|c|c|c}
				\hline
				$N$ & $y(x)$ in \eqref{eq:series1} & $y(x)$ in \eqref{eq:series2} & $y(x)$ in \eqref{eq:series1} & $y(x)$ in \eqref{eq:series2} \\
				& ($a=0.8$, $x=0.3$) & ($a=0.8$, $x=0.3$) & ($a=0.8$, $x=0.7$) & ($a=0.8$, $x=0.7$) \\
				\hline
				10   & 2.285559427400  & 2.276337892064  & 17.722665066666  & $1.00791\times 10^{5}$ \\
				50   & 2.285714285714  & 2.285714285695  & 26.622563574231  & $1.34009\times 10^{28}$ \\
				100  & 2.285714285714  & 2.285714285714  & 26.666611092450  & $1.99922\times 10^{57}$ \\
				200  & 2.285714285714  & 2.285714285714  & 26.666666666578  & $6.25120\times 10^{115}$ \\
				300  & 2.285714285714  & 2.285714285714  & 26.666666666667  & $2.25372\times 10^{174}$ \\
				400  & 2.285714285714  & 2.285714285714  & 26.666666666667  & $8.61497\times 10^{232}$ \\
				500  & 2.285714285714  & 2.285714285714  & 26.666666666667  & $3.40062\times 10^{291}$ \\
				600  & 2.285714285714  & 2.285714285714  & 26.666666666667  & $1.36992\times 10^{350}$ \\
				700  & 2.285714285714  & 2.285714285714  & 26.666666666667  & $5.59670\times 10^{408}$ \\
				800  & 2.285714285714  & 2.285714285714  & 26.666666666667  & $2.31012\times 10^{467}$ \\
				900  & 2.285714285714  & 2.285714285714  & 26.666666666667  & $9.61056\times 10^{525}$ \\
				1000 & 2.285714285714  & 2.285714285714  & 26.666666666667  & $4.02305\times 10^{584}$ \\
				\hline
			\end{tabular}
		}
	\end{table}
	
	For $x=0.7$, the table shows that $y(x)\approx 26.6667$ as $N\to\infty$, indicating convergence in Eq.~\eqref{eq:series1}, whereas Eq.~\eqref{eq:series2} demonstrates divergence. This discrepancy highlights the limitation of the absolute Convergence theorem for larger $x$ values.
	
\begin{figure}[h]
		\centering
		\includegraphics[scale=0.35]{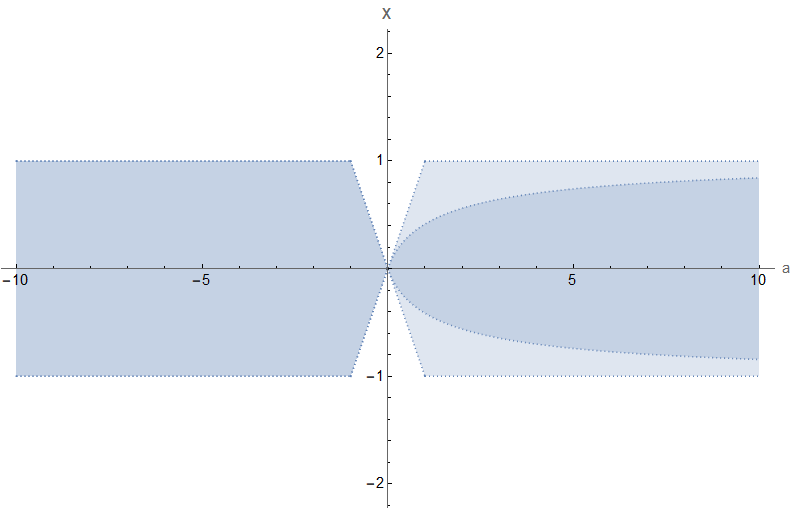}
		\caption{Domains of convergence for absolute and the radius of convergence under internal component rearrangement.}
		\label{choun3}
	\end{figure}

	In summary, Table~\ref{tab:merged} illustrates that the absolute Convergence theorem cannot reliably construct the radius of convergence for all cases, and the radius of convergence under internal component rearrangement  provides a more accurate description. Figure~\ref{choun3} visualizes the domains of convergence.
	\subsection{Absolute Convergence and Internal Component Rearrangement}\label{sec:convergence_analysis}
	
	A curious case arises as to why the P--P theorem cannot directly construct the radius of convergence for certain power series solutions. This section explores this issue by analyzing the difference between absolute convergence (based on the P--P theorem) and internal component rearrangement (achieved by rearranging the terms of the series).
	
	\paragraph{absolute convergence:} The P--P theorem is fundamentally built by observing the ratio $\overline{d}_{n+1}/\overline{d}_n$ in the recurrence relation as $n \to \infty$. This ratio corresponds to one of the roots of the characteristic polynomial of the recurrence. The radius of convergence for a power series $\sum_{n=0}^{\infty} d_n x^n$ is determined by requiring that the modulus of a root of the characteristic equation, multiplied by $|x|$, is less than unity:
	\begin{equation}
		\lim_{n\to\infty} \left|\frac{d_{n+1}}{d_n}\right| |x| < 1.
		\label{eq:conv_cond}
	\end{equation}
	
	For the local Heun function at $x=0$, the recurrence relation is given by:
	\begin{equation}
		\overline{d}_{n+1} = \alpha_1\overline{d}_n + \alpha_2\overline{d}_{n-1} = \frac{1+a}{a}\overline{d}_n - \frac{1}{a}\overline{d}_{n-1},
		\label{eq:heun_rec}
	\end{equation}
	where $\alpha_1 = (1+a)/a$ and $\alpha_2 = -1/a$. The closed-form solution for $\overline{d}_n$ is given by:
	\begin{equation}
		\overline{d}_n = \frac{-\left(\alpha_1-\sqrt{\alpha_1^2+4\alpha_2}\right)^{n+1} + \left(\alpha_1+\sqrt{\alpha_1^2+4\alpha_2}\right)^{n+1}}{2^{n+1}\sqrt{\alpha_1^2 + 4\alpha_2}}.
		\label{eq:closed_form}
	\end{equation}
	
	To evaluate the radius of convergence, multiply $\lim_{n\to\infty} |\overline{d}_{n+1}/\overline{d}_n|$ by $|x|$ and substitute \eqref{eq:closed_form}:
	\begin{equation}
		L = \lim_{n\to\infty} \left|\frac{\overline{d}_{n+1}}{\overline{d}_n}\right| |x| = \lim_{n\to\infty} \left|
		\frac{\frac{\alpha_1-\sqrt{\alpha_1^2+4\alpha_2}}{2} - \frac{\alpha_1+\sqrt{\alpha_1^2+4\alpha_2}}{2}\left(\frac{\alpha_1+\sqrt{\alpha_1^2+4\alpha_2}}{\alpha_1-\sqrt{\alpha_1^2+4\alpha_2}}\right)^{n+1}}{1 - \left(\frac{\alpha_1+\sqrt{\alpha_1^2+4\alpha_2}}{\alpha_1-\sqrt{\alpha_1^2+4\alpha_2}}\right)^{n+1}}
		\right| |x| < 1.
		\label{eq:radius_condition}
	\end{equation}
	
	Depending on the ratio of roots in \eqref{eq:radius_condition}, two cases arise:
	\begin{equation}
		\text{If}\ \left|\frac{\alpha_1+\sqrt{\alpha_1^2+4\alpha_2}}{\alpha_1-\sqrt{\alpha_1^2+4\alpha_2}}\right| < 1,\ \text{then}\ L = \frac{1}{2}\left|\alpha_1-\sqrt{\alpha_1^2+4\alpha_2}\right| |x| < 1,
		\label{eq:case1}
	\end{equation}
	\begin{equation}
		\text{If}\ \left|\frac{\alpha_1+\sqrt{\alpha_1^2+4\alpha_2}}{\alpha_1-\sqrt{\alpha_1^2+4\alpha_2}}\right| > 1,\ \text{then}\ L = \frac{1}{2}\left|\alpha_1+\sqrt{\alpha_1^2+4\alpha_2}\right| |x| < 1.
		\label{eq:case2}
	\end{equation}
	
	For the special case where $\alpha_1=2$ and $\alpha_2=-1$, \eqref{eq:radius_condition} simplifies to:
	\begin{equation}
		L = |x| < 1.
		\label{eq:simple_radius}
	\end{equation}
	\paragraph{Discrepancy with internal component rearrangement:}  
	By substituting $\alpha_1 = (1+a)/a$ and $\alpha_2 = -1/a$ into \eqref{eq:case1}--\eqref{eq:simple_radius}, we confirm that the radius of convergence for a local Heun function at $x=0$ matches the values in Table~\ref{cb.1}. However, the radius of convergence obtained when internal component rearrangement is applied, as shown in Table~\ref{cb.2}, is smaller. This discrepancy arises because internal component rearrangement allows the terms $|\overline{d}_n|$ to be split and reorganized, violating the structural invariance assumed in the P--P theorem.  
	
	\paragraph{Findings:}  
	While the P--P theorem provides a radius of convergence under the assumption of absolute convergence, it fails to account for the reduced radius observed when internal component rearrangement is applied. This highlights the necessity of imposing a no-rearrangement constraint to preserve the original radius of convergence, as demonstrated in this study.

	\subsection{Analysis of Hypergeometric Series and Radius of Convergence}\label{sec:hypergeometric_analysis}
	
	A hypergeometric function is defined by the power series:
	\begin{equation}
		\sum_{n=0}^{\infty} d_n x^n = 1 + \frac{ab}{c\,1!}x + \frac{a(a+1)b(b+1)}{c(c+1)\,2!}x^2 + \frac{a(a+1)(a+2)b(b+1)(b+2)}{c(c+1)(c+2)\,3!}x^3 + \cdots.
		\label{eq:hypergeom}
	\end{equation}
	
	This series consists of a two-term recurrence relation between successive coefficients. To determine the condition for absolute convergence of \eqref{eq:hypergeom}, we analyze the series of moduli $\sum_{n=0}^{\infty}|d_n||x|^n$:
	\begin{align}
		\sum_{n=0}^{\infty}|d_n||x|^n &= 1 + \left|\frac{ab}{c\,1!}\right||x| + \left|\frac{a(a+1)b(b+1)}{c(c+1)\,2!}\right||x|^2 \nonumber\\
		&\quad + \left|\frac{a(a+1)(a+2)b(b+1)(b+2)}{c(c+1)(c+2)\,3!}\right||x|^3 + \cdots.
		\label{eq:hypergeom_abs}
	\end{align}
	
	Applying the ratio test to \eqref{eq:hypergeom_abs}, the radius of convergence is obtained as:
	\begin{equation}
		\lim_{n\rightarrow \infty}\left|\frac{(n+a)(n+b)}{(n+c)(n+1)}\right| |x| < 1.
		\label{eq:ratio_test}
	\end{equation}
	
	Thus, $|x| < 1$ is required for the absolute convergence of the hypergeometric series.
	
	\paragraph{Asymptotic Expansion:} For the local Heun function, an asymptotic series expansion of \eqref{eq:heun_rec} is:
	\begin{align}
		\sum_{n=0}^{\infty} \overline{d}_n x^n &= 1 + \alpha_1x + (\alpha_1^2 + \alpha_2)x^2 + (\alpha_1^3 + 2\alpha_1\alpha_2)x^3 \nonumber\\
		&\quad + (\alpha_1^4 + 3\alpha_1^2\alpha_2 + \alpha_2^2)x^4 + (\alpha_1^5 + 4\alpha_1^3\alpha_2 + 3\alpha_1\alpha_2^2)x^5 + \cdots.
		\label{eq:asymp_exp}
	\end{align}
	\paragraph{Structural Constraints in the P--P Theorem:}
	The Poincar\'e--Perron (P--P) theorem traditionally provides a method for estimating the radius of convergence of a power series solution governed by a recurrence relation. However, this theorem assumes that the structural composition of recurrence terms remains invariant under any ordering.
	
	Our analysis demonstrates that in higher-order recurrence relations, the internal components of the recurrence terms can undergo rearrangement, leading to a reduction in the radius of convergence. We refer to this phenomenon as \textbf{Internal Component Rearrangement}, which violates the conventional assumptions of the P--P theorem.
	
	To address this issue, we impose a \textbf{Structural Constraint on the P--P Theorem}, explicitly prohibiting internal component rearrangement. This constraint ensures that the radius of convergence aligns with classical predictions, preserving the analytical validity of power series solutions. The necessity of this constraint highlights a fundamental limitation in the conventional application of the P--P theorem, reinforcing the role of structural coherence in determining convergence properties.
	
	The P--P theorem is typically applied to the ratio $\overline{d}_{n+1}/\overline{d}_n$ in the recurrence relation of a local Heun function at $x=0$. However, this approach is flawed when determining the radius of convergence. Instead, we must explicitly consider the absolute values of all terms in \eqref{eq:asymp_exp}:
	\begin{align}
		\sum_{n=0}^{\infty}|\overline{d}_n||x|^n &= 1 + |\alpha_1||x| + (|\alpha_1|^2 + |\alpha_2|)|x|^2 + (|\alpha_1|^3 + 2|\alpha_1||\alpha_2|)|x|^3 \nonumber\\
		&\quad + (|\alpha_1|^4 + 3|\alpha_1|^2|\alpha_2| + |\alpha_2|^2)|x|^4 \nonumber\\
		&\quad + (|\alpha_1|^5 + 4|\alpha_1|^3|\alpha_2| + 3|\alpha_1||\alpha_2|^2)|x|^5 + \cdots.
		\label{eq:abs_series}
	\end{align}
	
	To construct the radius of convergence, we replace $\alpha_1 \to |\alpha_1| = |(1+a)/a|$ and $\alpha_2 \to |\alpha_2| = |-1/a|$ in \eqref{eq:case1}--\eqref{eq:case2}. The radius of convergence for the local Heun function at $x=0$ then aligns with Table~\ref{cb.2}, except in the case of $a=-1$. This discrepancy arises because the P--P theorem assumes conditional convergence, whereas the absolute values of the coefficients in \eqref{eq:abs_series} must be considered for absolute convergence.
	
	\paragraph{Special Case $a = -1$:} For $a = -1$, substituting $\alpha_1 \to |\alpha_1| = |(1+a)/a|$ and $\alpha_2 \to |\alpha_2| = |-1/a|$ into \eqref{eq:radius_condition} yields:
	\begin{equation}
		\lim_{n \to \infty} \left|\frac{\frac{|A| - \sqrt{|A|^2 + 4|B|}}{2} - \frac{|A| + \sqrt{|A|^2 + 4|B|}}{2}\left(\frac{|A| + \sqrt{|A|^2 + 4|B|}}{|A| - \sqrt{|A|^2 + 4|B|}}\right)^{n+1}}{1 - \left(\frac{|A| + \sqrt{|A|^2 + 4|B|}}{|A| - \sqrt{|A|^2 + 4|B|}}\right)^{n+1}}\right| =  \lim_{n\rightarrow \infty } \left| \frac{-1+(-1)^n}{1+(-1)^n}\right|,
		\label{eq:special_case}
	\end{equation}
	which is undefined for determining convergence. Instead, by substituting $a = -1$ directly into \eqref{eq:17}, we can obtain the interval of convergence for the local Heun function at $x = 0$. 
	
	\paragraph{Findings:}  
	The Structural Constraints in the P--P Theorem and the Internal Component Rearrangement approach demonstrate that the radius of convergence depends critically on whether the structural composition of recurrence terms remains intact. While absolute convergence is preserved under conventional assumptions, it requires additional constraints to remain stable when internal component rearrangement is introduced, as shown in Table~\ref{cb.2}.  
	
	 \section{Domain Relationship between Absolute Convergence and Internal Component Rearrangement}  
	 \label{sec:domain_relation}
	 
	 According to Fig.~\ref{choun3}, the relationship between absolute convergence and internal component rearrangement is given by:  
	 \begin{scriptsize}
	 	\begin{equation}
	 		D_r = \textrm{Domain affected by internal component rearrangement} \subseteq D_o = \textrm{Domain of absolute convergence},
	 		\label{eq:domain_relation}
	 	\end{equation}
	 \end{scriptsize}
	 where \( D_o \) denotes the domain of convergence defined by the Poincaré–Perron theorem, and \( D_r \) represents the subdomain where internal component rearrangement alters the radius of convergence. This inclusion relationship holds not only for three-term recurrence relations in Fuchsian equations but also extends to general multi-term recurrence cases.
	 
	 \subsection{Uniqueness Theorem in the Framework of Fuchsian Equations}  
	 
	 The uniqueness theorem for differential equations asserts that a Fuchsian equation admits a unique solution under prescribed boundary conditions. For example, in a Poisson-type equation, if $\Omega$ is a bounded domain in $\mathbb{R}^n$, and a function $u: \Omega\rightarrow \mathbb{R}$ satisfies $u \in C^2(\Omega) \cap C(\bar{\Omega})$ and $\nabla^2 u = \rho$ in $\Omega$, along with either $u = f$ or $\partial u / \partial n = g$ on $\partial \Omega$, then $u$ is uniquely determined up to an additive con...
	 
	 A similar uniqueness principle applies to the Klein–Gordon equation, in which the scalar field $\Phi(r)$ defined in a region $\Omega$ is uniquely determined by Dirichlet or Neumann boundary conditions. This property naturally extends to the Heun equation, which arises in various physical contexts, such as the separation of variables in the $D=4$ Kerr–de Sitter metric~\cite{Cunh2015, Suzu1998, Tarl2016}. 
	 
	 In this setting, the uniqueness theorem ensures that the local Heun function is well-defined within a given analytic domain.
	 
	 \subsection{Uniqueness Breakdown Due to Internal Component Rearrangement}  
	 
	 Let us assume that the power series solution \( y^P(x) \) in \eqref{eq:ppsol}, derived from the Poincaré--Perron theorem at \( x=0 \), is unique within the domain \( D_o \), the domain of absolute convergence. Ideally, the alternative solution \( y^R(x) \), constructed via internal component rearrangement as in \eqref{eq:heun_series}, should coincide with \( y^P(x) \), possibly up to a constant multiple.
	 
	 However, in the region \( D_o - D_r \) (the bright shaded area in Fig.~\ref{choun3}, where \( a > 0 \)), the rearranged solution \( y^R(x) \) diverges. This divergence indicates that \( y^R(x) \) and \( y^P(x) \) become independent solutions, thereby violating the uniqueness theorem.
	 
	 This breakdown arises because the Poincaré--Perron theorem, when applied without structural constraints, fails to ensure invariance under internal component rearrangement. Thus, the uniqueness of solutions is not preserved in this framework unless absolute convergence is enforced from the outset.
	 
	 \subsection{Absolute Convergence as a Necessary Condition}  
	 
	 To restore the validity of the uniqueness theorem, let us regard \( y^R(x) \) as the uniquely defined solution within the structurally constrained domain \( D_r \). Since \( D_r \subseteq D_o \), the solution \( y^P(x) \) remains equivalent (or proportional) to \( y^R(x) \) within this subdomain.
	 
	 This relationship implies that \textbf{absolute convergence is a necessary condition} for the uniqueness of solutions to be maintained. By guaranteeing that all structurally admissible rearrangements converge to the same sum, absolute convergence preserves the functional identity of the solution series.
	 
	 Hence, imposing a \textbf{no-rearrangement constraint} reconciles the Poincaré--Perron theorem with the requirements of the uniqueness theorem. This synthesis aligns the mathematical rigor of analytic solutions with the physical expectation of structural consistency.

	\section{Uniqueness Theorem and the Necessity of Rearrangement Restriction}  
	\label{sec:uniqueness}  
	
	The uniqueness theorem is a cornerstone of mathematical physics, ensuring that solutions to differential equations, such as Poisson's equation, are uniquely determined by boundary conditions. However, the Poincaré--Perron (P--P) theorem, as traditionally formulated, does not inherently guarantee invariance under internal component rearrangement of series terms. This limitation can lead to a violation of the uniqueness condition. To resolve this inconsistency, we must explicitly impose a \textbf{no-rearrangement constraint} to ensure absolute convergence and preserve the structural integrity required for uniqueness.  
	
	\subsection{Illustration with Poisson's Equation}  
	
	Poisson's equation, which underlies Gauss's law in electrostatics, is expressed as:  
	\begin{equation}
		\mathbf{\nabla}\cdot(\epsilon_0\mathbf{\nabla}\varphi) = -\rho_f,
		\label{eq:poisson}
	\end{equation}
	where \( \varphi \) denotes the electric potential, and \( \mathbf{E} = -\mathbf{\nabla}\varphi \) is the corresponding electric field.  
	
	The uniqueness theorem asserts that for given boundary conditions, the solution \( \varphi \) is uniquely determined. If two potential solutions \( \varphi_1 \) and \( \varphi_2 \) exist, then their difference \( \phi = \varphi_2 - \varphi_1 \) satisfies Laplace's equation:  
	\begin{equation}
		\mathbf{\nabla}\cdot(\epsilon_0\mathbf{\nabla}\phi) = 0.
		\label{eq:laplace}
	\end{equation}
	Applying the divergence theorem yields \( \nabla\phi = 0 \), leading to \( \varphi_1 = \varphi_2 + C \), where \( C \) is a constant. Thus, the solution is unique up to an additive constant.  
	
	This fundamental principle extends to power series solutions, where absolute convergence ensures that internal rearrangement does not violate the uniqueness of solutions. The necessity of structural constraints, such as the no-rearrangement condition, becomes clear in this context.  
	
	\subsection{Connection Between Integrals and Series Convergence}  
	
\begin{figure}[!htb]
		\centering
		\minipage{0.4\textwidth}
		\includegraphics[width=\linewidth]{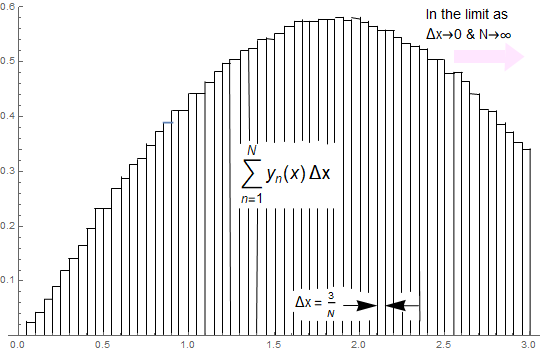}
		\caption{Partial sum}\label{integral1}
		\endminipage
		\hspace{0.02\textwidth}  
		\minipage{0.4\textwidth}
		\includegraphics[width=\linewidth]{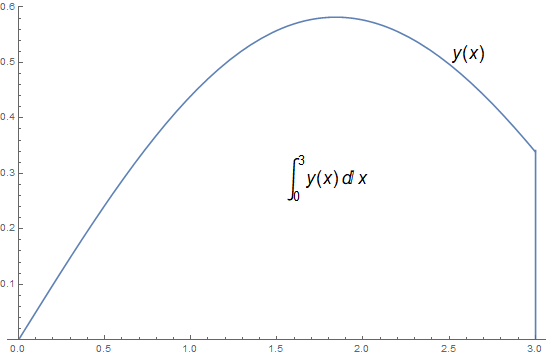}
		\caption{Integral}\label{integral2}
		\endminipage
	\end{figure}

	The uniqueness theorem is intimately linked to the invariance of integrals, which represent the limiting behavior of discrete sums. This relationship is depicted in Figs.~\ref{integral1} and \ref{integral2}.  
	
	As illustrated, partial sums approximate the integral, and as \( \Delta x \to 0 \), the sum converges to:  
	\begin{equation}
		\text{Area} = \int_{0}^{3} y(x)\,dx = \lim_{\Delta x \to 0} \sum_{n=1}^{N} y_n(x)\,\Delta x.
		\label{eq:integral}
	\end{equation}  
	
	This invariance mirrors the role of absolute convergence in power series: preserving structural integrity against term rearrangement.  
	
	\subsection{Invariance and Rearrangement Restriction}  
	
	The uniqueness theorem presumes that integrals remain invariant under rearrangement. However, the P--P theorem lacks this enforcement, making it vulnerable to breakdown under structural rearrangement.  
	
	To uphold uniqueness, we must distinguish \textbf{absolute convergence}, defined by P--P, from \textbf{internal component rearrangement}, which disrupts term structure. Only by imposing a no-rearrangement constraint can we ensure compatibility with the uniqueness theorem.  
	
	\subsection{Rearrangement and Uniqueness: Limitations of the Poincaré--Perron Theorem}  
	
	The P--P theorem assumes conditional convergence and fails to account for structural rearrangement. Even with the Laurent series expansion:  
	\begin{equation}
		f(z) = \sum_{n=-\infty}^{\infty} a_n(z - z_0)^n = \sum_{n=-\infty}^{\infty} b_n(z - z_0)^n,
		\label{eq:laurent}
	\end{equation}
	its uniqueness relies on:  
	\begin{equation}
		a_n = \frac{1}{2\pi i} \oint_C \frac{f(z)\,dz}{(z-z_0)^{n+1}}.
		\label{eq:cauchy}
	\end{equation}  
	
	Using:  
	\begin{equation}
		\oint_C (z - z_0)^{n-m-1}\,dz = 2\pi i\delta_{nm},
		\label{eq:residue}
	\end{equation}
	we conclude \( a_m = b_m \). But this conclusion presumes absolute convergence, which the P--P theorem does not guarantee.  
	
	\subsection{Necessity of Absolute Convergence for Integral and Series Exchange}  
	
	The interchange:  
	\begin{equation}
		\int_a^b \sum_{n=0}^\infty f_n\,dx \rightarrow \sum_{n=0}^\infty \int_a^b f_n\,dx,
		\label{eq:exchange}
	\end{equation}
	is only valid under absolute convergence. Otherwise, results may differ due to order-sensitive summation.  
	
	\subsection{Application to Heun Functions and Fuchsian Equations}  
	
	In 2004, Robert S. Maier reformulated the Heun equation \cite{Maie2007}:
	\begin{multline}
		Hl(a,q;\alpha,\beta,\gamma,\delta;x) = \left(1 - \frac{x}{a}\right)^{-\alpha} \\
		Hl\left(\frac{1}{1-a}, \frac{q - \gamma\alpha}{1-a}; -\beta + \gamma + \delta, \alpha, \gamma, \alpha - \beta + 1; \frac{x}{x-a}\right).
		\label{eq:heun_transform}
	\end{multline}
	
	Both sides must be absolutely convergent. Internal component rearrangement risks altering convergence behavior, again highlighting the need for structural constraints.  
	
	\subsection{Uniqueness Theorem and the Role of Absolute Convergence}  
	
	Absolute convergence maintains the validity of integration-series exchange and is vital for the uniqueness of solutions. Internal rearrangement disrupts this consistency. Therefore, a \textbf{no-rearrangement constraint} must be explicitly incorporated in the P--P framework.  
	
	\subsection{The Role of the Uniqueness Theorem and Frobenius Series}  
	
	In physical and mathematical modeling, discrete quantities approximate continuous structures. Although multiple series may emerge, the uniqueness theorem assures only one correct solution. The Frobenius series ensures this uniqueness by preserving absolute convergence.  
	
	Thus, absolute convergence is not a convenience—it is a necessity for consistency, accuracy, and structural coherence in differential equations.

	 \section{Conclusion and Discussion}  
	 
	 This paper has demonstrated that classical convergence theorems, such as Fuchs’ theorem and the Poincaré--Perron (P--P) theorem, implicitly rely on structural invariance in power series solutions. When this invariance is violated—specifically through internal component rearrangement of recurrence terms—convergence behavior changes fundamentally. To ensure uniqueness and consistency, we must impose a strict \textbf{no-rearrangement constraint}, preserving absolute convergence.
	 
	 \begin{cor}
	 	For a $d$-term recurrence relation derived from a Fuchsian differential equation with $d \geq 3$, the uniqueness theorem holds \emph{if and only if} the power series solution is absolutely convergent under a structural constraint. In the absence of such constraint, internal component rearrangement can violate uniqueness.
	 \end{cor}
	 
	 Consider the differential equation:
	 \begin{equation}
	 	y'' + p(z)y' + q(z)y = 0,
	 	\label{eq:diff_eq}
	 \end{equation}
	 with initial conditions \( y(0) = y_0 \), \( y'(0) = y_0' \). According to Fuchs' theorem, the radius of convergence of a Frobenius series solution is at least the minimum of the radii of convergence of \( p(z) \) and \( q(z) \). However, for a local Heun function, this prediction fails. As shown in Fig.~\ref{choun3}, the actual domain of convergence is often smaller.
	 
	 This discrepancy arises because the coefficients \( d_n \) are not atomic terms, but rather polynomials of degree \( n \) composed of nested recurrence components \( \alpha_{1,n}, \alpha_{2,n} \), etc. To evaluate convergence properly, the magnitude \( |d_n| \) must be treated as:
	 \begin{equation}
	 	|d_n| = \sum\prod_{i=0}^{j}\prod_{l=0}^{k}|\alpha_{1,i}||\alpha_{2,l}|.
	 	\label{eq:dn_abs}
	 \end{equation}
	 From this perspective, we conclude:
	 
	 \begin{cor}
	 	Fuchs' theorem, when applied to power series of the form \( y(z) = \sum_{n=0}^{\infty} d_n z^n \), remains valid only if each \( d_n \) is treated as an indivisible unit. For recurrence relations involving three or more terms, structural rearrangement leads to breakdowns in convergence and thus invalidates the theorem unless structural invariance is enforced.
	 \end{cor}
	 
	 \section{Further Discussion}  
	 
	 Lazarus Fuchs' theory, developed in the 19th century, remains foundational in the study of linear differential equations with regular singular points. Its elegance lies in transforming a differential equation into a recurrence relation by substituting a Frobenius-type power series \( y(z) = \sum_{n=0}^{\infty}d_n z^n \). Fuchs' theorem states that the radius of convergence of this series is at least as large as the minimum of the radii of convergence of the coefficients \( p(z) \) and \( q(z) \) in the differential equation. However, this classical result implicitly assumes that each coefficient \( d_n \) is treated as an indivisible entity—immune to structural rearrangement.
	 
	 In practice, particularly for multi-term recurrence relations as seen in the Heun equation, each \( d_n \) is not a single term but rather a composite of recurrence coefficients. If these internal components are rearranged or decomposed, the convergence properties change. Therefore, to ensure that the power series solution retains the convergence radius predicted by Fuchs' theorem, a strict \textbf{no-rearrangement constraint} must be imposed. This structural constraint is not an optional refinement; it is a necessary condition for preserving both analytic consistency and uniqueness.
	 
	 \vspace{0.3em}
	 \subsection*{A Geometric Analogy: Complex Numbers and Structural Integrity}
	 
	 To illustrate this principle, consider a familiar mathematical object: the complex number \( z = x + iy \). The power series \( \sum_{n=0}^\infty z^n \) converges absolutely for \( |z| < 1 \). Here, \( z \) is treated as a single geometric point in the complex plane. However, if one were to artificially split \( z \) into its real and imaginary parts and consider a constraint such as \( |x| + |y| < 1 \), the resulting convergence domain becomes strictly smaller than \( |z| < 1 \). This violates the geometric integrity of the complex number system.
	 
	 The same phenomenon occurs in recurrence relations. Just as we preserve the structure of \( z \) to ensure its algebraic and geometric properties, we must also preserve the ordering and indivisibility of recurrence components in \( d_n \). Otherwise, the convergence behavior—and thus the identity—of the function is compromised.
	 
	 \vspace{0.3em}
	 \subsection*{Analytical Integrity and Structural Order}
	 
	 This analogy reveals a broader structural insight. The convergence of a power series is not merely a matter of controlling the size of individual terms, but of maintaining the deep ordering embedded in the recurrence. Recurrence relations such as
	 \[
	 d_{n+1} = \alpha_{1,n} d_n + \alpha_{2,n} d_{n-1}
	 \]
	 contain within them a recursive memory. The coefficients \( d_n, d_{n-1}, \dots \) do not exist in isolation; their interplay encodes the differential equation’s structure. Disturbing this order—even slightly—alters the analytic landscape, much as rotating or stretching the complex plane would distort its geometry.
	 
	 In the case of Heun’s equation, which has four regular singularities and yields a three-term recurrence, the preservation of this structural integrity becomes essential. Rearrangement of internal components in such a setting does not merely affect numerical accuracy—it fundamentally breaks the conditions under which classical theorems like Fuchs' remain valid.
	 
	 \vspace{0.3em}
	 \subsection*{Toward a Unified Structural Principle}
	 
	 These observations lead to a more general conclusion: the radius of convergence is not merely a boundary on a complex plane—it is a \emph{structural invariant} of the system, reflecting deeper principles of symmetry, ordering, and indivisibility. Violating that structure by rearrangement undermines not only convergence but the uniqueness, consistency, and interpretability of solutions.
	 
	 Thus, we advocate the explicit introduction of a \textbf{no-rearrangement constraint} as a foundational axiom when extending Fuchsian theory to modern differential equations, particularly those involving complex recurrence structures such as the Heun equation. This constraint ensures not only the mathematical rigor of convergence estimates but also their compatibility with physical and geometric intuition.
	 
	 In this light, absolute convergence is more than a technical condition—it is a guarantor of structural coherence. Just as complex numbers derive their power from their inseparability, so too do power series solutions derive their stability from an indivisible internal structure. The no-rearrangement constraint is, therefore, a reflection of this deeper mathematical truth.

	 \section{Findings and Potential Applications}  
	 
	 In this study, we investigated the radii of convergence for power series solutions to Fuchsian differential equations, focusing on how they are affected by internal component rearrangement of the recurrence structure. Specifically, we contrasted the classical notion of \textbf{absolute convergence}—as derived via Fuchs’ theorem—with the reduced radius obtained when internal rearrangement is permitted. Our findings reveal that the radius of convergence under such rearrangement is not determined by the explicit coefficients of the equation, but rather by the relative geometric configuration of the singularities.
	 
	 These observations are not only mathematically significant but also conceptually illuminating. They underscore the necessity of a \textbf{no-rearrangement constraint} to preserve structural integrity, ensure uniqueness, and maintain compatibility with classical convergence theorems. This insight opens new directions for research and applications in mathematical analysis and theoretical physics.
	 
	 \subsection*{Key Implications}
	 
	 \begin{enumerate}
	 	\item \textbf{Structural Stability of Power Series Solutions:} Imposing a no-rearrangement constraint protects the convergence behavior of solutions by preventing the distortion of recurrence structure, thereby upholding the uniqueness theorem.
	 	
	 	\item \textbf{Refined View of Singularities and Analytic Domains:} The results suggest that the influence of singular points on the domain of convergence is more subtle than classical theory predicts, especially for higher-order recurrence relations.
	 	
	 	\item \textbf{Generalization to Higher-Order Systems:} These findings naturally extend to differential equations with more than three terms, where recurrence structures become more intricate and sensitive to structural alterations.
	 \end{enumerate}
	 
	 \subsection*{Future Research Directions}
	 
	 Building upon these findings, several promising avenues of research emerge:
	 
	 \begin{itemize}
	 	\item \textbf{Extension to Complex Differential Systems:} Explore how structural constraints affect convergence in systems with multiple or irregular singularities.
	 	
	 	\item \textbf{Asymptotics and Stokes Phenomena:} Investigate the connection between internal rearrangement and asymptotic behavior, particularly how it influences Stokes transitions and analytic continuation.
	 	
	 	\item \textbf{Partial Differential Equations:} Extend the framework to PDEs, focusing on the role of structural invariance in separable or eigenfunction-based solutions.
	 	
	 	\item \textbf{Analytic-Symbolic Duality:} Develop formal methods to classify when a given recurrence relation necessitates structural preservation, and identify scenarios where symbolic manipulations can misrepresent the true analytic domain.
	 \end{itemize}
	 
	 This study establishes a foundation for a broader structural framework in the analysis of differential equations—one that integrates geometric intuition, recurrence theory, and analytic rigor. The no-rearrangement constraint, far from being a technical limitation, emerges as a unifying principle connecting convergence, uniqueness, and the intrinsic architecture of analytic series.

	\appendix
	\section{Convergence Domains Under Internal Component Rearrangement for 192 Local Heun Functions}  
	
	A systematically generated classification of 192 local solutions of the Heun equation, which exhibit an isomorphic structure to the Coxeter group associated with the Coxeter diagram \( D_4 \), was obtained by Maier (2007) \cite{Maie2007}.  
	
	We explicitly compute the modified convergence domains due to internal component rearrangement for nine representative cases among the 192 local solutions of the Heun equation, as presented in Table 2 \cite{Maie2007}.

	\subsection{ ${\displaystyle (1-x)^{1-\delta } Hl(a, q - (\delta  - 1)\gamma a; \alpha - \delta  + 1, \beta - \delta + 1, \gamma ,2 - \delta ; x)}$ \\  \emph{and}\;  ${\displaystyle x^{1-\gamma } (1-x)^{1-\delta } Hl(a, q-(\gamma +\delta -2)a-(\gamma -1)(\alpha +\beta -\gamma -\delta +1); \alpha - \gamma -\delta +2}$ \\ ${\displaystyle, \beta - \gamma -\delta +2, 2-\gamma, 2 - \delta ; x)}$}
	boundary conditions of the independent variable $x$  of  $Hl(a, q - (\delta  - 1)\gamma a; \alpha - \delta  + 1, \beta - \delta + 1, \gamma ,2 - \delta ; x)$ and
	$Hl(a, q-(\gamma +\delta -2)a-(\gamma -1)(\alpha +\beta -\gamma -\delta +1); \alpha - \gamma -\delta +2, \beta - \gamma -\delta +2, 2-\gamma, 2 - \delta ; x)$ are the same as  Eq. (\ref{eq:17}).
	\subsection{ ${\displaystyle  Hl(1-a,-q+\alpha \beta; \alpha,\beta, \delta, \gamma; 1-x)}$ \\
		\emph{and}  ${\displaystyle (1-x)^{1-\delta } Hl(1-a,-q+(\delta -1)\gamma a+(\alpha -\delta +1)(\beta -\delta +1); \alpha-\delta +1,\beta-\delta +1}$\\
		${\displaystyle, 2-\delta, \gamma; 1-x)}$}
	Replace coefficient $a$ and independent variable $x$ by $1-a$ and $1-x$, respectively, in (\ref{eq:17}).
	
	The convergence condition under internal component rearrangement is given by
	\begin{equation}
		\left|\frac{1}{1-a}(1-x)^2\right| +\left|\frac{2-a}{1-a}(1-x)\right|<1 \hspace{.5cm}\mbox{where}\;\; a \ne 1 \nonumber 
	\end{equation}
	\subsection{ ${\displaystyle x^{-\alpha } Hl\left(\frac{1}{a},\frac{q+\alpha [(\alpha -\gamma -\delta +1)a-\beta +\delta ]}{a}; \alpha , \alpha -\gamma +1, \alpha -\beta +1,\delta ;\frac{1}{x}\right)}$}
	Replace coefficient $a$ and $x$ by $\frac{1}{a}$ and $\frac{1}{x}$, respectively, in (\ref{eq:17}).
	The convergence condition under internal component rearrangement is given by
	\begin{equation}
		\left| a x^{-2}\right| +\left| (1+a)x^{-1} \right|<1 \nonumber
	\end{equation}
	\subsection{ ${\displaystyle \left(1-\frac{x}{a} \right)^{-\beta } Hl\left(1-a, -q+\gamma \beta; -\alpha +\gamma +\delta, \beta, \gamma, \delta; \frac{(1-a)x}{x-a} \right)}$ \\
		\emph{and} \small ${\displaystyle (1-x)^{1-\delta }\left(1-\frac{x}{a} \right)^{-\beta+\delta -1} Hl\left(1-a, -q+\gamma [(\delta -1)a+\beta -\delta +1]; -\alpha +\gamma +1, \beta -\delta+1, \gamma, 2-\delta; \frac{(1-a)x}{x-a} \right)}$\normalsize}
	Replace $a$ and $x$ by $1-a$ and ${\displaystyle \frac{(1-a)x}{x-a}}$, respectively, in (\ref{eq:17}).
	The Condition of Convergence Under Internal Component Rearrangementof is given by
	\begin{equation}
		\left|\frac{(1-a)x^2}{(x-a)^2}\right|+\left|\frac{(2-a)x}{(x-a)}\right|<1 \hspace{.5cm}\mbox{where}\;\; x \ne a \nonumber 
	\end{equation}
	\subsection{ ${\displaystyle x^{-\alpha } Hl\left(\frac{a-1}{a}, \frac{-q+\alpha (\delta a+\beta -\delta )}{a}; \alpha, \alpha -\gamma +1, \delta , \alpha -\beta +1; \frac{x-1}{x} \right)}$}
	Replace $a$ and $x$ by $\frac{a-1}{a}$ and $\frac{x-1}{x}$, respectively, in (\ref{eq:17}).
	The convergence condition under internal component rearrangement is given by
	\begin{equation}
		\left|\frac{a}{(1-a)}\frac{(x-1)^2}{x^2} \right| +\left|\frac{(1-2a)}{(1-a)}\frac{(x-1)}{x} \right|<1 \hspace{.5cm}\mbox{where}\;\; a \ne 1\nonumber
	\end{equation}
	\subsection{ ${\displaystyle \left(\frac{x-a}{1-a} \right)^{-\alpha } Hl\left(a, q-(\beta -\delta )\alpha ; \alpha , -\beta+\gamma +\delta , \delta , \gamma; \frac{a(x-1)}{x-a} \right)}$}
	Replace $x$ by $\frac{a(x-1)}{x-a}$ in (\ref{eq:17}).
	The convergence condition under internal component rearrangement is given by
	\begin{equation}
		\left| \frac{a(x-1)^2}{(x-a)^2}\right| +\left|\frac{(1+a)(x-1)}{(x-a)} \right|<1 \hspace{.5cm}\mbox{where}\;\; x \ne a\nonumber
	\end{equation}
	\section{Domains of Resummed Convergence of 4 Fuchsian differential equations having multi-term recurrence relations}
	\subsection{Riemann's P--differential equation}
	The Riemann P--differential equation is written by
	\begin{eqnarray}
		&&\frac{d^2{y}}{d{z}^2} + \left( \frac{1-\alpha -\alpha^{'}}{z-a} +\frac{1-\beta -\beta^{'}}{z-b} +\frac{1-\gamma -\gamma^{'}}{z-c} \right) \frac{d{y}}{d{z}} \label{eq:73}\\
		&&+  \left( \frac{\alpha \alpha^{'}(a-b)(a-c)}{z-a} +\frac{\beta\beta^{'}(b-c)(b-a)}{z-b} +\frac{\gamma \gamma^{'}(c-a)(c-b)}{z-c} \right)\frac{y}{(z-a)(z-b)(z-c)} = 0 \nonumber
	\end{eqnarray}
	With the condition $\alpha+ \alpha^{'}+ \beta+ \beta^{'}+ \gamma+ \gamma^{'}= 1$. The regular singular points are $a$, $b$ and $c$ with exponents $\{ \alpha, \alpha^{'} \}$, $\{ \beta, \beta^{'} \}$ and $\{ \gamma, \gamma^{'} \}$. This equation was first obtained in the form by Papperitz. \cite{Barn1908,Papp1885}
	The solutions to the Riemann P--differential equation are given in terms of the hypergeometric function.   However, currently, the analytic solutions of Riemann's differential equation at $a$, $b$ and $c$ are unknown because of its complicated mathematical calculation: the five-term recurrence relation of a Frobenius series starts to appear. 
	
	We assume the solution takes the form
	\begin{equation}
		y(x)= \sum_{n=0}^{\infty } c_n x^{n+\lambda }  \label{eq:75}
	\end{equation}
	By putting (\ref{eq:75}) into (\ref{eq:73}) at $x=z-a$, the algebraic equation around $z=a$ is
	\begin{equation}
		\frac{d^2{y}}{d{x}^2} + \left( \frac{1-\overline{\alpha }}{x} +\frac{1-\overline{\beta }}{x+k} +\frac{1-\overline{\gamma }}{x+l} \right) \frac{d{y}}{d{x}}
		+  \left( \frac{\alpha \alpha^{'}kl}{x} -\frac{\beta\beta^{'}mk}{x+k} +\frac{\gamma \gamma^{'}ml}{x+l} \right)\frac{y}{x(x+k)(x+l)} = 0 \label{eq:76}
	\end{equation}
	For convenience for calculations, the parameters have been changed by setting
	\begin{equation}
		k=a-b,\;\; l=a-c,\;\; m=b-c,\;\; \overline{\alpha }=\alpha +\alpha ^{'},\;\; \overline{\beta }=\beta +\beta ^{'},\;\; \overline{\gamma }=\gamma +\gamma ^{'} \label{eq:77}
	\end{equation}
	with $\overline{\alpha }+\overline{\beta }+\overline{\gamma }=1$.
	We obtain the recurrence system by substituting (\ref{eq:75}) into (\ref{eq:76}).
	\begin{equation}
		c_{n+1}=\alpha _{1,n} \;c_n +\alpha _{2,n} \;c_{n-1} +\alpha _{3,n} \;c_{n-2} +\alpha _{4,n} \;c_{n-3}\hspace{1cm};n\geq3
		\label{eq:78}
	\end{equation}
	where,
	\begin{subequations}
		\begin{equation}
			\alpha _{1,n} =  -\frac{(n+\lambda )\left( 2(k+l)(n+\lambda )-k(\overline{\alpha }-\overline{\beta })-l(\overline{\alpha }-\overline{\gamma })\right) +(k+l)\alpha \alpha ^{'}-m(\beta\beta ^{'}-\gamma \gamma ^{'})}{kl(n+1-\alpha +\lambda )(n+1-\alpha^{'} +\lambda)}
			\label{eq:79a}
		\end{equation}
		\begin{equation}
			\alpha _{2,n} = - \frac{(n-1+\lambda )\left( (k^2+4kl+l^2)(n-1+\lambda )+k^2\overline{\beta }+l^2 \overline{\gamma }+2kl(1-\overline{\alpha })\right) +kl\alpha \alpha^{'}-kl\beta \beta^{'}+ml\gamma \gamma ^{'}}{k^2 l^2(n+1-\alpha +\lambda )(n+1-\alpha^{'} +\lambda)}
			\label{eq:79b}
		\end{equation}
		\begin{equation}
			\alpha _{3,n} = - \frac{(n-2+\lambda )\left( 2(k+l)(n-2+\lambda) +k(1+\overline{\beta })+ l(1+\overline{\gamma })\right)}{k^2 l^2(n+1-\alpha +\lambda )(n+1-\alpha^{'} +\lambda)}
			\label{eq:79c}
		\end{equation}
		\begin{equation}
			\alpha _{4,n} = - \frac{(n-3+\lambda )(n-2+\lambda )}{k^2 l^2(n+1-\alpha +\lambda )(n+1-\alpha^{'} +\lambda)}
			\label{eq:79d}
		\end{equation}
		and
		\begin{equation}
			c_1= \alpha _{1,0} c_0,\;\; c_2= \left( \alpha _{1,0}\alpha _{1,1} +\alpha _{2,1}\right) c_0,\;\; c_3= \left( \alpha _{1,0}\alpha _{1,1}\alpha _{1,2} +\alpha _{1,0}\alpha _{2,2}+ \alpha _{1,2}\alpha _{2,1}+ \alpha _{3,2}\right) c_0
			\label{eq:79e}
		\end{equation}
	\end{subequations}
	We have two indicial roots which are $\lambda =\alpha $ and $\alpha ^{'}$. (\ref{eq:78}) is a 5-term recurrence relation as we know.
	
	The convergence condition under internal component rearrangement  of (\ref{eq:75}) is obtained by applying  Thm.\ref{main} such as
	\begin{equation}
		\left|\frac{2(2a-b-c)}{(a-b)(a-c)}x \right| + \left|\frac{(2a-b-c)^2+2(a-b)(a-c)}{(a-b)^2 (a-c)^2} x^2 \right| + \left|\frac{2(2a-b-c)}{(a-b)^2 (a-c)^2} x^3 \right| + \left|\frac{1}{(a-b)^2 (a-c)^2} x^4 \right| <1  \nonumber
	\end{equation}
	\subsection{Heine differential equation}
	The Heine differential equation is defined by
	\begin{equation}
		\frac{d^2{y}}{d{z}^2} + \frac{1}{2}\left( \frac{1}{z-a_1} +\frac{2}{z-a_2} +\frac{2}{z-a_3} \right) \frac{d{y}}{d{z}}
		+ \frac{1}{4}\left( \frac{\beta _0+ \beta _1 z+ \beta _2 z^2+ \beta _3 z^3 }{(z-a_1)(z-a_2)^2 (z-a_3)^2} \right) y= 0 \label{eq:97}
	\end{equation}
	It has four regular singular points at $a_1$, $a_2$, $a_3$ and $\infty $. Parameters $a_2$ and $a_3$ are identical to each other \cite{Moon1961,Zwil1997}. 
	The differential equation around $z=a_1$, which is obtained from (\ref{eq:97}) by setting $x=z-a_1$ is
	\begin{equation}
		\frac{d^2{y}}{d{x}^2} + \frac{1}{2}\left( \frac{1}{x} +\frac{2}{x-m} +\frac{2}{x-k} \right) \frac{d{y}}{d{x}}
		+ \frac{\beta _3 x^3 + \overline{\beta _2} x^2 + \overline{\beta _1} x + \overline{\beta _0}}{4x(x-m)^2(x-k)^2}y = 0 \label{eq:98}
	\end{equation}
	For convenience for calculations, the parameters have been changed by setting
	\begin{equation}
		m=a_2-a_1,\;\; k=a_3-a_1,\;\; \overline{\beta _2} =3a_1 \beta_3 + \beta _2,\;\; \overline{\beta _1}= 3a_1^2 \beta _3 + 2a_1 \beta _2 +\beta _1,\;\; \overline{\beta _0} = a_1^3 \beta _3 + a_1^2 \beta _2 +a_1\beta _1 + \beta _0 \nonumber 
	\end{equation}
	Looking for a solution of (\ref{eq:98}) through the expansion
	\begin{equation}
		y(x) = \sum_{n=0}^{\infty } c_n x^{n+\lambda }\nonumber 
	\end{equation}
	We obtain the following conditions:
	\begin{equation}
		c_{n+1}=\alpha _{1,n} \;c_n +\alpha _{2,n} \;c_{n-1} +\alpha _{3,n} \;c_{n-2} +\alpha _{4,n} \;c_{n-3}\hspace{1cm};n\geq3
		\label{kk:100}
	\end{equation}
	where,
	\begin{subequations}
		\begin{equation}
			\alpha _{1,n} = \frac{2mk(m+k)(n+\lambda )^2 -\frac{\overline{\beta _0}}{4}}{m^2 k^2(n+1 +\lambda )(n+\frac{1}{2} +\lambda)}
			\label{eq:101a}
		\end{equation}
		\begin{equation}
			\alpha _{2,n} = - \frac{(m^2+4mk+k^2)(n-1+\lambda )(n-\frac{1}{2}+\lambda )+ \frac{\overline{\beta _1}}{4}}{m^2 k^2(n+1 +\lambda )(n+\frac{1}{2} +\lambda)}
			\label{eq:101b}
		\end{equation}
		\begin{equation}
			\alpha _{3,n} = \frac{2(m+k)(n-2+\lambda )(n-1+\lambda )+ \frac{\overline{\beta _2}}{4}}{m^2 k^2(n+1+\lambda )(n+\frac{1}{2} +\lambda)}
			\label{eq:101c}
		\end{equation}
		\begin{equation}
			\alpha _{4,n} = -\frac{(n-3+\lambda )(n-\frac{3}{2}+\lambda )+\frac{\overline{\beta _3}}{4}}{m^2 k^2(n+1 +\lambda )(n+\frac{1}{2} +\lambda)}
			\label{eq:101d}
		\end{equation}
		and
		\begin{equation}
			c_1= \alpha _{1,0} c_0,\;\; c_2= \left( \alpha _{1,0}\alpha _{1,1} +\alpha _{2,1}\right) c_0,\;\; c_3= \left( \alpha _{1,0}\alpha _{1,1}\alpha _{1,2} +\alpha _{1,0}\alpha _{2,2}+ \alpha _{1,2}\alpha _{2,1}+ \alpha _{3,2}\right) c_0
			\label{eq:101e}
		\end{equation}
	\end{subequations}
	We have two indicial roots which are $\lambda =0$ and $\frac{1}{2}$. (\ref{kk:100}) is the 5-term recurrence relation.
	
	The convergence condition under internal component rearrangement  of  the Heine equation at $z= a_1$ is obtained by applying  Thm.\ref{main} such as
	\begin{eqnarray}
		&&\left|\frac{2(2a_1- a_2 -a_3)}{(a_2-a_1)(a_3-a_1)}x  \right| + \left|\frac{(2a_1- a_2 - a_3)^2+ 2(a_2-a_1)(a_3-a_1)}{(a_2-a_1)^2 (a_3-a_1)^2}x^2 \right| +\left| \frac{2(2a_1- a_2 - a_3)}{(a_2-a_1)^2 (a_3-a_1)^2}x^3 \right| \nonumber\\
		&&+ \left|\frac{1}{(a_2-a_1)^2 (a_3-a_1)^2}x^4 \right|<1  \nonumber
	\end{eqnarray}
	The differential equation around $z=a_2$, which is obtained from (\ref{eq:97}) by setting $x=z-a_2$. Putting a power series $y(x)= \sum_{n=0}^{\infty }c_n x^{n+\lambda }$ into the new (\ref{eq:97}), we obtain the 4-term recurrence relation as follows.
	\begin{equation}
		c_{n+1}=\alpha _{1,n} \;c_n +\alpha _{2,n} \;c_{n-1} +\alpha _{3,n} \;c_{n-2} \hspace{1cm};n\geq2
		\label{eq:119}
	\end{equation}
	where,
	\begin{subequations}
		\begin{equation}
			\alpha _{1,n} = -\frac{k(2m+k)(n+\lambda )(n+\frac{1}{2}+\lambda ) +\frac{\overline{\beta _1}}{4}}{m k^2(n+1 +\lambda )^2 + \frac{\overline{\beta _0}}{4}}
			\label{eq:120a}
		\end{equation}
		\begin{equation}
			\alpha _{2,n} = - \frac{(m+2 k)(n-1+\lambda )(n+\lambda ) +\frac{\overline{\beta _2}}{4}}{m k^2(n+1+\lambda )^2 + \frac{\overline{\beta _0}}{4}}
			\label{eq:120b}
		\end{equation}
		\begin{equation}
			\alpha _{3,n} = - \frac{(n-2+\lambda )(n-\frac{1}{2}+\lambda )+ \frac{\overline{\beta _3}}{4}}{m k^2(n+1+\lambda )^2 + \frac{\overline{\beta _0}}{4}}
			\label{eq:120c}
		\end{equation}
		and
		\begin{equation}
			c_1= \alpha _{1,0} c_0,\;\; c_2= \left( \alpha _{1,0}\alpha _{1,1} +\alpha _{2,1}\right) c_0
			\label{eq:120d}
		\end{equation}
	\end{subequations}
	where
	\begin{equation}
		m=a_2-a_1,\;\; k=a_2-a_3,\;\; \overline{\beta _2} =3a_2 \beta_3 + \beta _2,\;\; \overline{\beta _1}= 3a_2^2 \beta _3 + 2a_2 \beta _2 +\beta _1,\;\; \overline{\beta _0} = a_2^3 \beta _3 + a_2^2 \beta _2 +a_2\beta _1 + \beta _0 \nonumber 
	\end{equation}
	We have two indicial roots which are $\lambda =\pm\sqrt{\frac{-\overline{\beta _0}}{4mk^2}}$.
	
	The convergence condition under internal component rearrangement  of the Heine   equation at $z=a_2$ is obtained by applying  Thm.\ref{main} such as
	\begin{equation}
		\left| \frac{(a_2-a_3)(3a_2-2a_1-a_3)}{(a_2-a_1)(a_2-a_3)^2}x \right| +\left| \frac{(3a_2- a_1 - 2a_3)}{(a_2-a_1)  (a_2-a_3)^2}x^2\right| + \left| \frac{1}{(a_2-a_1)  (a_2-a_3)^2}x^3 \right|<1  \nonumber
	\end{equation}
	\subsection{Generalized Heun's differential equation}
	The generalized Heun equation (GHE) is a second-order linear ordinary differential equation of the form 
	\begin{equation}
		\frac{d^2{y}}{d{x}^2} + \left( \frac{1-\mu _0}{x} +\frac{1-\mu _1}{x-1} +\frac{1-\mu _2}{x-a} +\alpha \right) \frac{d{y}}{d{x}}
		+ \frac{\beta _0+ \beta _1 x+ \beta _2 x^2 }{x(x-1)(x-a)} y= 0 \label{eq:130}
	\end{equation}
	where $\mu _0, \mu _1, \mu _2, \alpha, \beta _0, \beta _1, \beta _2 \in \mathbb{C}$. It has three regular singular points which are 0, 1 and $a$ with exponents $\{0,\mu _0\}$, $\{0,\mu _1\}$ and $\{0,\mu _2\}$, while $\infty $ is at most an irregular singularity. This equation was first introduced by Sch$\ddot{\mbox{a}}$fke and Schmidt. \cite{Schaf1980}
	Heun's equation is derived from the generalized Heun's differential equation by changing all coefficients $\alpha =\beta _2 =0$, $1-\mu _0 =\gamma $, $1-\mu _1 =\delta $, $1-\mu _2 =\epsilon $, $\beta _1 =\alpha \beta $ and $\beta _0= -q$. 
	
	Looking for a solution of (\ref{eq:130}) through the expansion
	\begin{equation}
		y(x) = \sum_{n=0}^{\infty } c_n x^{n+\lambda }\nonumber 
	\end{equation}
	We obtain the following conditions:
	\begin{equation}
		c_{n+1}=\alpha _{1,n} \;c_n +\alpha _{2,n} \;c_{n-1} +\alpha _{3,n} \;c_{n-2}\hspace{1cm};n\geq 2
		\label{eq:131}
	\end{equation}
	where,
	\begin{subequations}
		\begin{equation}
			\alpha _{1,n} = \frac{(n+\lambda )\left( (1+a)(n+\lambda )+1-\mu _0-\mu _2 +a(1-\mu _0-\mu _1 -\alpha )\right) -\beta _0}{a(n+1 +\lambda )(n+1-\mu _0 +\lambda)}
			\label{eq:132a}
		\end{equation}
		\begin{equation}
			\alpha _{2,n} = - \frac{(n-1+\lambda )\left( n+\lambda +1-\mu _0-\mu _1-\mu _2 -(1+a)\alpha \right) +\beta _1}{a(n+1 +\lambda )(n+1-\mu _0 +\lambda)}
			\label{eq:132b}
		\end{equation}
		\begin{equation}
			\alpha _{3,n} = - \frac{\alpha (n-2+\lambda ) +\beta _2}{a(n+1 +\lambda )(n+1-\mu _0 +\lambda)}
			\label{eq:132c}
		\end{equation}
		and
		\begin{equation}
			c_1= \alpha _{1,0} c_0,\;\; c_2= \left( \alpha _{1,0}\alpha _{1,1} +\alpha _{2,1}\right) c_0
			\label{eq:132d}
		\end{equation}
	\end{subequations}
	We have two indicial roots which are $\lambda =0$ and $\mu _0$. (\ref{eq:131}) is the 4-term recurrence relation.
	
	The convergence condition under internal component rearrangement  of  the GHE is obtained by applying  Thm.\ref{main} such as 
	\begin{equation}
		\left| \frac{1+a}{a} x\right| + \left|- \frac{1}{a} x^2\right|<1  \label{zzz}
	\end{equation}

	The GHE around $x=1$ is taken by putting $z=1-x$ into (\ref{eq:130}).
	\begin{equation}
		\frac{d^2{y}}{d{z}^2} + \left( \frac{1-\mu _1}{z} +\frac{1-\mu _0}{z-1} +\frac{1-\mu _2}{z-(1-a)} -\alpha \right) \frac{d{y}}{d{z}}
		+ \frac{-(\beta _0 +\beta _1 +\beta _2 )+ (\beta _1 +2\beta _2)z - \beta _2 z^2 }{z(z-1)(z-(1-a))} y= 0 \label{eq:142}
	\end{equation}
	If we compare (\ref{eq:142}) with (\ref{eq:130}), all coefficients on the above are correspondent to the following way.
	\begin{equation}
		\begin{split}
			& a \longrightarrow   1-a \\ & \mu _0 \longrightarrow  \mu _1 \\ & \mu _1 \longrightarrow  \mu _0 \\
			& \alpha  \longrightarrow  -\alpha \\ & \beta _0 \longrightarrow -(\beta _0 +\beta _1 +\beta _2) \\
			& \beta _1 \longrightarrow \beta _1 +2\beta _2 \\ & \beta _2 \longrightarrow -\beta _2 \\ & x \longrightarrow z = 1-x
		\end{split}\label{eq:143}
	\end{equation}
	Putting (\ref{eq:143}) into  (\ref{zzz}), we obtain domain of convergence Under internal component rearrangement of the GHE around $x=1$ where $z= 1-x$ for an infinite series.
	
	The GHE around $x=a$ is taken by putting $z=\frac{a-x}{a}$ into (\ref{eq:130}).
	\begin{equation}
		\frac{d^2{y}}{d{z}^2} + \left( \frac{1-\mu _2}{z} +\frac{1-\mu _0}{z-1} +\frac{1-\mu _1}{z-(a^{-1}-1)} -a\alpha \right) \frac{d{y}}{d{z}}
		+ \frac{-(\beta _1 +a\beta _2 +a^{-1}\beta _0 )+ (\beta _1 +2 a\beta _2)z - a\beta _2 z^2 }{z(z-1)(z-(a^{-1}-1))} y= 0 \label{eq:144}
	\end{equation}
	If we compare (\ref{eq:144}) with (\ref{eq:130}), all coefficients on the above are correspondent to the following way.
	\begin{equation}
		\begin{split}
			& a \longrightarrow   a^{-1}-1 \\ & \mu _0 \longrightarrow  \mu _2 \\ & \mu _1 \longrightarrow  \mu _0 \\ & \mu _2 \longrightarrow  \mu _1 \\
			& \alpha  \longrightarrow  -a\alpha \\ & \beta _0 \longrightarrow -(\beta _1 +a\beta _2 +a^{-1}\beta _0) \\
			& \beta _1 \longrightarrow \beta _1 +2 a\beta _2 \\ & \beta _2 \longrightarrow -a\beta _2 \\ & x \longrightarrow z = \frac{a-x}{a}
		\end{split}\label{eq:145}
	\end{equation}
	Putting (\ref{eq:145}) into  (\ref{zzz}), we obtain domain of convergence Under internal component rearrangement of the GHE around $x=a$ where $z= \frac{a-x}{a}$ for an infinite series.
	\subsection{The Lam\'{e} (or ellipsoidal) wave equation}
	The Lam\'{e} (or ellipsoidal) wave equation arises from deriving Helmholtz equation in ellipsoidal coordinates $\nabla^2 \phi + \omega \; \phi =0$ where $\nabla^2$ is the Laplacian, $\omega $ is the wavenumber and $\phi $ is the amplitude. \cite{ Arsc1983,Arsc1956,Arsc1959,Erde1955} Whereas Lam\'{e} equation is derived from separation of the Laplace equation in confocal ellipsoidal coordinates $\nabla^2 \phi =0$, only the special case of the ellipsoidal wave equation.
	
	The ellipsoidal wave equation is a second-order linear ODE of the Jacobian form 
	\begin{equation}
		\frac{d^2{y}}{d{z}^2} - \left( a+b k^2 sn^2z +q k^4 sn^4z \right) y= 0 \label{eq:146}
	\end{equation}
	where the Jacobian elliptic function $sn\;z = sn(z,k)$ has modulus $k$. If $q=0$, it turns to be the Lame differential equation.
	If we take $x= sn^2 z$ as an independent variable in (\ref{eq:146}), we obtain the algebraic form of it.
	\begin{equation}
		\frac{d^2{y}}{d{x}^2} + \frac{1}{2}\left( \frac{1}{x} +\frac{1}{x-1} +\frac{1}{x-\delta } \right) \frac{d{y}}{d{x}}
		+ \frac{ \beta + \mu x +\gamma  x^2 }{x(x-1)(x-\delta )} y= 0 \label{eq:147}
	\end{equation}
	Here, the parameters have been changed by setting 
	\begin{equation}
		\delta =k^{-2}, \;\;\beta = -\frac{a}{4k^2}, \;\;\mu =-\frac{b}{4}, \;\;\gamma =-\frac{k^2}{4}q \nonumber
	\end{equation}
	
	Looking for a solution of (\ref{eq:147}) through the expansion
	\begin{equation}
		y(x) = \sum_{n=0}^{\infty } c_n x^{n+\lambda }\nonumber 
	\end{equation}
	We obtain the following conditions:
	\begin{equation}
		c_{n+1}=\alpha _{1,n} \;c_n +\alpha _{2,n} \;c_{n-1} +\alpha _{3,n} \;c_{n-2}\hspace{1cm};n\geq 2
		\label{eq:148}
	\end{equation}
	where,
	\begin{subequations}
		\begin{equation}
			\alpha _{1,n} = \frac{(1+\delta )(n+\lambda )^2 -\beta }{\delta (n+1 +\lambda )\left( n+\frac{1}{2} +\lambda \right)}
			\label{eq:149a}
		\end{equation}
		\begin{equation}
			\alpha _{2,n} = - \frac{(n-1+\lambda )\left( n-\frac{1}{2}+\lambda \right) +\mu }{\delta (n+1 +\lambda )\left( n+\frac{1}{2} +\lambda \right)}
			\label{eq:149b}
		\end{equation}
		\begin{equation}
			\alpha _{3,n} = - \frac{\gamma }{\delta (n+1 +\lambda )\left( n+\frac{1}{2} +\lambda \right)}
			\label{eq:149c}
		\end{equation}
		and
		\begin{equation}
			c_1= \alpha _{1,0} c_0,\;\; c_2= \left( \alpha _{1,0}\alpha _{1,1} +\alpha _{2,1}\right) c_0
			\label{eq:149d}
		\end{equation}
	\end{subequations}
	We have two indicial roots which are $\lambda =0$ and $\frac{1}{2}$. (\ref{eq:148}) is the 4-term recurrence relation. 
	
	The convergence condition under internal component rearrangement of  the Lam\'{e} wave equation around $x=0$ is obtained by applying  Thm.\ref{main} such as 
	\begin{equation}
		\left| \frac{1+\delta}{\delta} x\right| +\left| \frac{-1}{\delta} x^2\right|<1 \label{zzz.1}
	\end{equation} 
	The Lam\'{e} wave equation around $x=1$ is taken by putting $\varrho =1-x$ into (\ref{eq:147}).
	\begin{equation}
		\frac{d^2{y}}{d{\varrho}^2} + \frac{1}{2}\left( \frac{1}{\varrho} +\frac{1}{\varrho -1} +\frac{1}{\varrho -(1-\delta )} \right) \frac{d{y}}{d{\varrho}}
		+ \frac{-(\beta + \gamma ) + (\mu +2\gamma )\varrho  - \gamma \varrho ^2 }{\varrho (\varrho -1)(\varrho -(1-\delta ))} y= 0 \label{eq:159}
	\end{equation}
	If we compare (\ref{eq:159}) with (\ref{eq:147}), all coefficients on the above are correspondent to the following way.
	\begin{equation}
		\begin{split}
			& \delta  \longrightarrow   1-\delta  \\ & \beta  \longrightarrow -(\beta +\gamma ) \\ & \mu \longrightarrow  \mu +2\gamma \\
			& \gamma   \longrightarrow  -\gamma \\ & x \longrightarrow z = 1-\varrho
		\end{split}\label{eq:160}
	\end{equation}
	Putting (\ref{eq:160}) into (\ref{zzz.1}), we obtain domain of convergence Under internal component rearrangement of the Lam\'{e} wave equation around $x=1$ where $\varrho= 1-x$.
	
	The Lam\'{e} wave equation around $x=\delta $ is taken by putting $\phi =1-\delta ^{-1}x$ into (\ref{eq:147}).
	\begin{equation}
		\frac{d^2{y}}{d{\phi}^2} + \frac{1}{2}\left( \frac{1}{\phi} +\frac{1}{\phi -1} +\frac{1}{\phi -(1-\delta ^{-1} )} \right) \frac{d{y}}{d{\phi}}
		+ \frac{-( \gamma \delta +\mu +\delta ^{-1}\beta ) + (\mu +2\gamma \delta )\phi - \gamma \delta \phi ^2 }{\phi (\phi -1)(\phi -(1-\delta ))} y= 0 \label{eq:161}
	\end{equation}
	If we compare (\ref{eq:161}) with (\ref{eq:147}), all coefficients on the above are correspondent to the following way.
	\begin{equation}
		\begin{split}
			& \delta  \longrightarrow   1-\delta ^{-1}  \\ & \beta  \longrightarrow -( \gamma \delta +\mu +\delta ^{-1}\beta) \\ & \mu \longrightarrow  \mu +2\gamma \delta \\
			& \gamma   \longrightarrow  -\gamma \delta \\ & x \longrightarrow \phi = 1-\delta ^{-1}x
		\end{split}\label{eq:162}
	\end{equation}
	Putting (\ref{eq:162}) into (\ref{zzz.1}), we obtain domain of convergence Under internal component rearrangement of the Lam\'{e} wave equation around $x=\delta $ where $\phi = 1-\delta ^{-1} x$. 
	
\end{document}